\documentclass[11pt]{article}
\usepackage{mathtools}
\usepackage{amsmath,amsthm}
\usepackage{csquotes}

\usepackage{amssymb}

\usepackage{graphicx}

\usepackage{subfig}

\usepackage{mathpazo}

\usepackage{soul}

\usepackage{xcolor}

\usepackage{cancel}

\usepackage{empheq}

\usepackage{lipsum}

\definecolor{myblue}{rgb}{.8, .8, 1}
  \newcommand*\mybluebox[1]{
    \colorbox{myblue}{\hspace{1em}#1\hspace{1em}}}

\usepackage[obeyspaces,hyphens,spaces]{url}

\usepackage{hyperref}
\hypersetup{
    colorlinks=true,
    linkcolor=blue,
    citecolor=blue,
    filecolor=magenta,
    urlcolor=cyan
}

\usepackage[
  open,
  openlevel=2,
  atend,
  numbered
]{bookmark}

\usepackage[capitalize, nameinlink, noabbrev]{cleveref}
\crefname{equation}{}{}
\crefname{chapter}{Chapter}{Chapters}
\crefname{item}{item}{items}
\crefname{figure}{Figure}{Figures}
\crefname{theorem}{Theorem}{Theorems}
\crefname{lemma}{Lemma}{Lemmas}
\crefname{proposition}{Proposition}{Propositions}
\crefname{corollary}{Corollary}{Corollarys}
\crefname{definition}{Definition}{Definitions}
\crefname{fact}{Fact}{Facts}
\crefname{example}{Example}{Examples}
\crefname{algorithm}{Algorithm}{Algorithms}
\crefname{remark}{Remark}{Remarks}
\crefname{note}{Note}{Notes}
\crefname{notation}{Notation}{Notations}
\crefname{case}{Case}{Cases}
\crefname{exercise}{Exercise}{Exercises}
\crefname{question}{Question}{Questions}
\crefname{claim}{Claim}{Claims}
\crefname{enumi}{}{}

\usepackage[top= 2cm, bottom = 2 cm, left = 2.2 cm, right= 2.2 cm]{geometry}

\usepackage{float}

\numberwithin{equation}{section}

\theoremstyle{plain}
\newtheorem{theorem}{Theorem}[section]

\newtheorem{corollary}[theorem]{Corollary}
\newtheorem{fact}[theorem]{Fact}
\newtheorem{lemma}[theorem]{Lemma}
\newtheorem{proposition}[theorem]{Proposition}

\theoremstyle{definition}
\newtheorem{definition}[theorem]{Definition}
\newtheorem{example}[theorem]{Example}

\newtheorem{remark}[theorem]{Remark}

\usepackage{enumitem}

\newcommand{\inte}{\ensuremath{\operatorname{int}}}

\newcommand{\cone}{\ensuremath{\operatorname{cone}}}

\newcommand{\spn}{\ensuremath{{\operatorname{span} \,}}}
\newcommand{\weakly}{\ensuremath{{\;\operatorname{\rightharpoonup}\;}}}

\newcommand{\dist}{\ensuremath{\operatorname{d}}}

\newcommand{\Pro}{\ensuremath{\operatorname{P}}}

\providecommand{\abs}[1]{\lvert#1\rvert}
\providecommand{\norm}[1]{\lVert#1\rVert}
\providecommand{\Norm}[1]{{\Big\lVert}#1{\Big\rVert}}

\providecommand{\innp}[1]{\langle#1\rangle}
\providecommand{\Innp}[1]{\Big\langle#1\Big\rangle}

\begin{document}

\title{ \sffamily  On  angles between convex sets  in Hilbert spaces}

\author{
         Heinz H.\ Bauschke\thanks{
                 Mathematics, University of British Columbia, Kelowna, B.C.\ V1V~1V7, Canada.
                 E-mail: \href{mailto:heinz.bauschke@ubc.ca}{\texttt{heinz.bauschke@ubc.ca}}.},~
         Hui\ Ouyang\thanks{
                 Mathematics, University of British Columbia, Kelowna, B.C.\ V1V~1V7, Canada.
                 E-mail: \href{mailto:hui.ouyang@alumni.ubc.ca}{\texttt{hui.ouyang@alumni.ubc.ca}}.},~
         and Xianfu\ Wang\thanks{
                 Mathematics, University of British Columbia, Kelowna, B.C.\ V1V~1V7, Canada.
                 E-mail: \href{mailto:shawn.wang@ubc.ca}{\texttt{shawn.wang@ubc.ca}}.}
                 }

\date{April 3, 2021}

\maketitle

\begin{abstract}
The notion of the angle between two subspaces has a long history, dating 
back to Friedrichs's work in 1937 and Dixmier's work on the minimal angle in 1949. 
In 2006, Deutsch and Hundal studied extensions to convex sets in order 
to analyze convergence rates for the cyclic projections algorithm. 

In this work, we characterize the positivity of the minimal angle between two
convex cones. We show the existence
of, and necessary conditions for, optimal solutions of minimal angle problems associated with two convex subsets as well. Moreover, we generalize a result by Deutsch on minimal angles  from  linear subspaces to cones. This generalization yields sufficient conditions for the closedness of the sum of two closed convex cones. This also relates to conditions proposed by
Beutner and by Seeger and Sossa.
Furthermore, we investigate the relation between the
intersection of two cones (at least one of which is nonlinear) and the intersection of the polar and dual cones of the underlying cones. It turns out that the two angles involved 
cannot be  positive simultaneously. 
Various examples illustrate the sharpness of our results. 
\end{abstract}

{\small

}

\section{Introduction} \label{sec:Introduction}

Throughout this paper, we shall assume that
\begin{empheq}[box = \mybluebox]{equation*}
\text{$\mathcal{H}$ is a real Hilbert space},
\end{empheq}
with inner product $\innp{\cdot,\cdot}$ and induced norm $\|\cdot\|$.
Moreover, we assume that $\mathcal{H} \neq {0}$ and that $m \in \mathbb{N}
\smallsetminus \{0\}$, where $\mathbb{N}=\{0,1,2,\ldots\}$.

The angle and the minimal angle between two closed linear subspaces were first
introduced by Friedrichs \cite{Friedrichs1937} in 1937 and by Dixmier
\cite{Dixmier1949} in 1949, respectively. (We recommend the nice expository
paper \cite{Deutsch1995} for details on the angle and minimal angle between
linear subspaces.) In order to characterize the rate of convergence of the
cyclic projections algorithm for the intersection of finitely many closed
convex sets in Hilbert spaces, Deutsch and Hundal generalized these
definitions from linear subspaces to general convex sets in \cite{DH2006II}.

Motivated by the applications of the angle and minimal angle between two
convex cones presented in \cite{DH2006II}, \cite{LewisLukeMalick2009}, \cite{Obert1991},  \cite{SS2016I}, and \cite{Tenenhaus1988} (see the last section for more details on applications), we show the existence of, and necessary conditions for, 
optimal solutions of minimal angle problems (the optimal solution is called a \emph{pair of principal vectors} in this work), 
and provide results on the evaluation of the cosine of the minimal angle.
Let us present our main results: 
\begin{enumerate}
	\item[\textbf{R1:}]  
	In \cref{theorem:c0}, we provide equivalent expressions for the cosine of
	the minimal angle between two closed convex sets. This yields 
	characterizations of the positivity of the minimal angle 
	(see \cref{theorem:c0K1K2EQ}).

	\item[\textbf{R2:}] 	
As a generalization of \cite[Theorem~2.12(1)$\Rightarrow$(2)]{Deutsch1995},
\cref{theorem:K1K2closed} states that if the minimal angle between two nonempty closed
convex cones is positive, then the vector
difference of the cones is closed. As an application of
\cref{theorem:K1K2closed}, we provide in \cref{cor:K1K2closed} sufficient
conditions for the closedness of the sum of two closed convex cones. In particular, some of
these conditions reduce to the sufficient conditions in
\cite[Theorem~3.2]{Beutner2007} and \cite[Proposition~4.1]{SS2016I}.

	\item[\textbf{R3:}] 
			\cref{theorem:KominusMperpNeq0} illustrates that, for nonempty closed convex cones $K_{1}$ and $K_{2}$ with $K_{1}$ nonlinear,
			under some
			assumptions (e.g., when $\mathcal{H}$ is finite-dimensional),
			$K_{1} \cap K_{2} =\{0\}$ implies that $K^{\ominus}_{1} \cap
			K^{\oplus}_{2} \neq \{0\}$, which shows that
			$c_{0}(K_{1},K_{2})<1$ and $c_{0}(K_{1}^{\ominus},
			K_{2}^{\oplus}) <1$ cannot occur simultaneously.
\end{enumerate}

The paper is organized as follows.   In  \cref{sec:Preliminaries}, we present some basic results that will be used in the next sections. In \cref{sec:Principalvectors}, we  provide equivalent expressions for the cosine of the minimal angle between two nonempty convex sets and show
 properties of the pair of principal vectors. 
Using the equivalent expressions for the cosine of the minimal angle between two nonempty convex sets given in \cref{sec:Principalvectors}, we characterize the  cosine of the minimal angle between two nonempty closed convex cones being strictly less than $1$ in  \cref{sec:anglesclosedconvexcones}. Moreover, 
in \cref{sec:anglesclosedconvexcones},  we  also give sufficient conditions for the closedness of the sum of two closed convex cones and consider the relation between the intersection of nonempty closed convex cones $K_{1}$ and $K_{2}$ and the intersection of $K_{1}^{\ominus}$ and $K_{2}^{\oplus}$.   Our applications and conclusions are stated in \cref{sec:applications}.

We now turn to the notation used in this work. Set  $\mathbf{B}_{\mathcal{H}}:=\{x\in \mathcal{H} ~:~ \norm{x} \leq 1 \}$ and  $\mathbf{S}_{\mathcal{H}}:=\{x \in \mathcal{H} ~:~ \norm{x}=1 \}$. Denote by  $\mathbb{R}_{+}:=\{\lambda \in \mathbb{R} ~:~ \lambda \geq 0 \}$ and $\mathbb{R}_{++}:=\{\lambda \in \mathbb{R} ~:~ \lambda >0 \}$.  Let $D$ be a nonempty subset
of $\mathcal{H}$. 
$\overline{D}$ is the \emph{closure} of $D$.
$D$ is a \emph{cone} if $D=\mathbb{R}_{++}D$. 
The \emph{conical hull} of $D$ is the intersection of all the cones in $\mathcal{H}$ containing $D$, i.e., the smallest cone in $\mathcal{H}$ containing $D$. It is denoted by $\cone D$. The \emph{closed conical hull} of $D$ is the smallest closed cone in $\mathcal{H}$ containing $D$.  It is denoted by $\overline{\cone} D$. Because we shall use results on conical hull in both \cite{BC2017} and \cite{D2012}, we point out that  by \cite[Proposition~6.2]{BC2017}, when $D$ is a convex set, we have that the definition  of conical hull in  \cite[Definition~6.1]{BC2017} presented above and the one shown in \cite[page~45]{D2012} are consistent, and that  $\overline{\cone} D= \overline{\cone D}$.
The \emph{polar cone} (or \emph{negative dual cone}) of $D$ is the set $D^{\ominus} := \{ u \in \mathcal{H} ~:~ (\forall x \in D) \innp{x,u} \leq 0  \}$.
The \emph{(positive) dual cone} and \emph{orthogonal complement} of $D$ are   $D^{\oplus} :=-D^{\ominus}$ and $ D^{\perp} :=D^{\ominus}  \cap D^{\oplus} =\{u \in \mathcal{H}~:~ (\forall x \in D) \innp{x,u}=0 \} $, respectively.
Let $C$ be a nonempty closed convex subset of $\mathcal{H}$. The \emph{projector} (or \emph{projection operator}) onto $C$ is the operator, denoted by
$\Pro_{C}$,  that maps every point in $\mathcal{H}$ to its unique projection onto $C$. 
Moreover, $(\forall x \in \mathcal{H})$ $\dist_{C}  (x) :=\min_{c \in C} \norm{x-c} =\norm{x -\Pro_{C}x}$. 
Suppose that $D$ is a nonempty closed convex subset of $\mathcal{H}$ as well. Denote $\dist (C, D) := \inf \{ \norm{x-y} ~:~ x \in C, y \in D \} $.   A sequence $(x_{k})_{k \in \mathbb{N}}$ in $\mathcal{H}$ \emph{converges
	weakly} to a point $x \in \mathcal{H}$ if, for every $u \in \mathcal{H}$,
$\innp{x_{k},u} \rightarrow \innp{x,u}$; in symbols, $x_{k} \weakly x$.

For other notation not explicitly defined here, we refer the reader to \cite{BC2017}.

\section{Preliminaries} \label{sec:Preliminaries}
In this section, we collect some results that will be used subsequently. 

\subsection*{Cones and projectors}
\begin{fact} {\rm \cite[Proposition~6.2(i)]{BC2017} } \label{fact:coneC}
	Let $C$ be a  subset of  $\mathcal{H}$. Then $\cone C = \mathbb{R}_{++} C$.
\end{fact}

\begin{fact} {\rm \cite[Propositions~6.3 and 6.4]{BC2017} } \label{lemma:K}
	Let $K$ be a nonempty convex cone.  Then  $K+K=K$. Moreover, if $-K \subseteq K$, then $K$ is a   linear subspace.
\end{fact}

\begin{fact} {\rm \cite[Theorem~4.5]{D2012}} \label{fact:dualcone}
	Let $C$ be a nonempty subset of  $\mathcal{H}$.  Then the following hold:
	\begin{enumerate}
		\item  \label{fact:dualcone:o} $C^{\ominus}$ is a closed convex cone and $C^{\perp}$ is a closed linear subspace. 
		\item   \label{fact:dualcone:eq}$C^{\ominus} =(\overline{C})^{\ominus}=\left( \cone(C) \right)^{\ominus}=\left( \overline{\cone}(C) \right)^{\ominus}$.
		\item \label{fact:dualcone:oo} $ C^{\ominus\ominus} = \overline{\cone}(C) $.

		\item  \label{fact:dualcone:closedconvex} If $C$ is a closed convex cone, then $C^{\ominus\ominus} =C$.
		\item  \label{fact:dualcone:linearsubsp} Assume that $C$ is a linear subspace. Then $C^{\ominus}=C^{\perp}$. In addition, if $C$ is closed, then $C=C^{\ominus\ominus}=C^{\perp \perp}$.
	\end{enumerate}
\end{fact}

 \begin{fact} {\rm \cite[Theorem~4.6]{D2012}} \label{fact:K1:Km:dualsum}
 	Let $K_{1}, \ldots, K_{m}$ be closed convex cones in $\mathcal{H}$. Then $\left( \cap^{m}_{i=1} K_{i} \right)^{\ominus}=\overline{\sum^{m}_{i=1} K^{\ominus}_{i}}$.
 \end{fact}

\begin{lemma} \label{lemma:polardual}
	Let $C$  be a nonempty subset  of $\mathcal{H}$.   Then the following statements hold:
	\begin{enumerate}
		\item \label{lemma:polardual:-C} $	(-C)^{\ominus} =- C^{\ominus} =C^{\oplus}$.
		\item \label{lemma:polardual:oplus} $C^{\oplus\oplus}=\overline{\cone} (C)$.
	
		\item  \label{lemma:polardual:llinear} Assume that $C$ is a linear subspace of $\mathcal{H}$. Then $C^{\perp}=C^{\ominus}=C^{\oplus}$.
	\end{enumerate}	
\end{lemma}
\begin{proof}
	\cref{lemma:polardual:-C}: By definition of polar cone and dual cone, $- C^{\ominus} =C^{\oplus}$. Let $x \in \mathcal{H}$. Then 
	\begin{align*}
	x \in (-C)^{\ominus}  \Leftrightarrow (\forall c \in C) \quad \innp{x, -c} \leq 0 \Leftrightarrow (\forall c \in C) \quad \innp{-x, c} \leq 0  \Leftrightarrow -x \in C^{\ominus}  \Leftrightarrow x \in -C^{\ominus}, 
	\end{align*}
	which implies that $(-C)^{\ominus} =- C^{\ominus}$.
	
	\cref{lemma:polardual:oplus}: Using \cref{lemma:polardual:-C} above and \cref{fact:dualcone}\cref{fact:dualcone:oo}, we see that  $C^{\oplus \oplus} =( C^{\oplus})^{\oplus} =(-C^{\ominus})^{\oplus} =C ^{ \ominus \ominus}=\overline{\cone} (C)$.
	
	\cref{lemma:polardual:llinear}: By \cref{fact:dualcone}\cref{fact:dualcone:o}$\&$\cref{fact:dualcone:linearsubsp},  clearly
	$C^{\oplus} = -C^{\ominus}=-C^{\perp} = C^{\perp} =C^{\ominus}$.
\end{proof}

\begin{fact} {\rm \cite[Theorem~3.16]{BC2017}} \label{fact:characterPC}
Let $C$ be a nonempty closed convex subset of $\mathcal{H}$. Let $x \in \mathcal{H}$. Then there is exactly one best approximation to $x$  from $C$.   Moreover,  for every $p$ in  $ \mathcal{H}$,
	\begin{align*}
	p=\Pro_{C} x \Leftrightarrow \left[ p \in C ~ \text{ and }~ (\forall y \in C) ~ \innp{y-p, x-p} \leq 0\right].
	\end{align*}
\end{fact}

\begin{fact}  {\rm \cite[Propositions~6.28]{BC2017}} \label{fact:chara:PK}
	Let $K$ be a nonempty closed convex cone in $\mathcal{H}$, let $x \in \mathcal{H}$, and let
 $p \in \mathcal{H}$. Then $ p =\Pro_{K} x \Leftrightarrow \left[  p \in K, x-p \perp p,  ~ \text{and}~ x-p \in K^{\ominus}\right]$.
\end{fact}

\subsection*{Angles between convex sets}

\begin{definition} {\rm \cite[Definitions~2.3 and 3.2]{DH2006II}} \label{defn:Angles}
Let $C$ and $D$ be   nonempty convex sets in $\mathcal{H}$. The \emph{minimal angle} between $C$ and $D$ is the angle in $\left[0, \frac{\pi}{2} \right]$ whose cosine is given by 
\begin{align} \label{eq:minimalangle}
c_{0} (C,D) := \sup \left\{   \innp{x,y} ~:~ x \in \overline{\cone} (C) \cap \mathbf{B}_{\mathcal{H}}, ~ y \in \overline{\cone} (D) \cap \mathbf{B}_{\mathcal{H}}\right\}.
\end{align}
In addition, the \emph{angle} between $C$ and $D$ is the angle in $\left[0, \frac{\pi}{2} \right]$ whose cosine is given by 
\begin{align*}
c (C,D) := c_{0} \big( (\cone C) \cap (\overline{ C^{\ominus} +D^{\ominus}}), ~ (\cone D) \cap  (\overline{ C^{\ominus} +D^{\ominus}}) \big).
\end{align*}
\end{definition}

Because the cosine function is decreasing, the angle corresponding to the optimal solution of \cref{eq:minimalangle} is the minimal one in question. To facilitate statements, we refer to the corresponding constrained maximization problem to \cref{eq:minimalangle} as the \emph{minimal angle problem}.
 
\begin{fact} {\rm  \cite[page~48]{SS2016I} } \label{fact:coBS}
	Let $C$ and $D$ be two nonempty convex sets of $\mathcal{H}$ such that $C \neq \{0\}$ and $D \neq \{0\}$.  Then 
	\begin{align*}
	c_{0} (C,D) = \max \left\{ 0, \sup \left\{   \innp{x,y} ~:~ x \in \overline{\cone} (C) \cap \mathbf{S}_{\mathcal{H}}, ~ y \in \overline{\cone} (D) \cap \mathbf{S}_{\mathcal{H}}\right\}   \right\}.
	\end{align*}
\end{fact}

According to \cref{fact:coBS}, the minimal angle  defined in \cref{defn:Angles}
 is different from the \enquote{minimal angle} considered in  \cite{Obert1991}, \cite{SS2016I}, and \cite{Tenenhaus1988}, but  when both  angles are in $\left[0, \frac{\pi}{2}\right]$, they are consistent.  Because we care mainly on using the cosine of minimal angle to describe  convergence rates of algorithms, and the convergence rate is always nonnegative, in this work we only study the minimal angle defined in \cref{defn:Angles}.

Because   for every 
nonempty convex set $C$, $\overline{\cone} (C)=\overline{\cone} (\overline{C})=\overline{\cone} (\cone C) $, although \cite[Lemma~2.4]{DH2006II} shows only  $ c_{0} (C,D)=c_{0} \left( \overline{\cone} (C) ,\overline{\cone} (D) \right)$, in the following \cref{fact:AnglesProperties}\cref{fact:AnglesProperties:c0:EQ} we present $  c_{0} (C,D)=c_{0} \left( \overline{ C} ,\overline{D} \right)=c_{0} \left(  \cone  (C) , \cone (D) \right)=c_{0} \left( \overline{\cone} (C) ,\overline{\cone} (D) \right)$.  Therefore, it is clear that when $\overline{ C^{\ominus} +D^{\ominus}} =\mathcal{H}$, then $c (C,D) =c_{0} \left(  \cone  (C) , \cone (D) \right)= c_{0}(C,D)  $.

\begin{fact} {\rm \cite[Lemma~2.4, Theorem~2.5 and Proposition~3.3]{DH2006II}}
\label{fact:AnglesProperties}
	Let $C$ and $D$ be   nonempty convex sets of $\mathcal{H}$. Then:
	\begin{enumerate}
		\item \label{fact:AnglesProperties:between0and1}$ c_{0} (C,D)  \in \left[0,1\right]$ and $ c (C,D)  \in \left[0,1\right]$. 
		\item \label{fact:AnglesProperties:ineq} $(\forall x \in \overline{\cone} (C) )$ $(\forall y \in \overline{\cone} (D) )$ $\innp{x,y} \leq   c_{0} (C,D) \norm{x} \norm{y}$.
		\item  \label{fact:AnglesProperties:c0:EQ} $ c_{0} (C,D)=c_{0} (D, C)=c_{0} \left( \overline{ C} ,\overline{D} \right)=c_{0} \left(  \cone  (C) , \cone (D) \right)=c_{0} \left( \overline{\cone} (C) ,\overline{\cone} (D) \right)$.
	\end{enumerate}
\end{fact}

\begin{lemma} \label{lemma:CD}
		Let $C$ and $D$ be   nonempty convex subsets of $\mathcal{H}$. Then the following hold:
		\begin{enumerate}
			\item \label{lemma:CD:UV} Let $U$ and $V$ be   nonempty convex subsets of $\mathcal{H}$ such that $C \subseteq U$ and $D \subseteq V$. Then $c_{0}(C,D) \leq c_{0}(U,V)$. 
			\item \label{lemma:CD:--}  $c_{0}(C,D) = c_{0}(-C,-D)$, $c_{0}(-C,D) = c_{0}(C,-D)$, $c (C,D) = c (-C,-D)$, and  $c (-C,D) = c (C,-D)$. 

			\item  \label{lemma:CD:C0=1} Suppose that $ \big(\overline{\cone} (C) \cap \overline{\cone} (D) \big) \smallsetminus \{0\} \neq \varnothing$. Then $c_{0} (C,D)  =1$.

			\item \label{lemma:CD:c0:c} $0 \leq c(C,D) \leq c_{0} (C,D) \leq 1$.
		\end{enumerate}
\end{lemma}

\begin{proof}
	\cref{lemma:CD:UV}: This follows easily from \cref{defn:Angles}.
	
	\cref{lemma:CD:--}: Note that for every nonempty convex subset $A$ of $\mathcal{H}$, $\overline{\cone} (-A) =- \overline{\cone} A$. Hence, by \cref{defn:Angles}, it is easy to see that 
	\begin{align} \label{eq:lemma:CD:--}
	c_{0}(C,D) = c_{0}(-C,-D)\quad \text{and} \quad c_{0}(-C,D) = c_{0}(C,-D).
	\end{align}
	In addition, by \cref{lemma:polardual}\cref{lemma:polardual:-C}, clearly,
	$ \overline{ (-C)^{\ominus} +(-D)^{\ominus}} =\overline{ - (C^{\ominus} +D^{\ominus})} =-\left(\overline{ C^{\ominus} +D^{\ominus}}\right)$. Hence, by \cref{eq:lemma:CD:--} and \cref{defn:Angles}, we obtain that $c (C,D) = c (-C,-D)$, and  $c (-C,D) = c (C,-D)$. 
	
		\cref{lemma:CD:C0=1}: Take $x \in \left(  \overline{\cone} (C) \cap \overline{\cone} (D) \right) \smallsetminus \{0\} $.  Then by \cref{fact:coneC}, $   \frac{x}{\norm{x}} \in \left( \overline{\cone} (C) \cap \mathbf{B}_{\mathcal{H}}  \right) \cap \left(\overline{\cone} (D)  \cap \mathbf{B}_{\mathcal{H}} \right) $. Combine this with \cref{defn:Angles} and  \cref{fact:AnglesProperties}\cref{fact:AnglesProperties:between0and1} to obtain that $	1 \geq c_{0} (C,D) \geq \Innp{ \frac{x}{\norm{x}} , \frac{x}{\norm{x}}  }=1$,
	which implies that $c_{0} (C,D)  =1$.

	 \cref{lemma:CD:c0:c}: By \cref{fact:AnglesProperties}\cref{fact:AnglesProperties:c0:EQ}  and 	\cref{lemma:CD:UV} above, 
	 \begin{align*}
c_{0}(C,D) =c_{0} \left(  \cone  (C) , \cone (D) \right) \geq c_{0} \left( (\cone C) \cap (\overline{ C^{\ominus} +D^{\ominus}}) , ~ (\cone D) \cap  (\overline{ C^{\ominus} +D^{\ominus}}) \right) =c (C,D),
	 \end{align*}
	 which, combining with \cref{fact:AnglesProperties}\cref{fact:AnglesProperties:between0and1}, deduces  the desired results.
\end{proof}

\begin{remark} \label{remark:c}
	Below, we provide three particular examples of  four-tuples $(C,D,U,V)$ of closed convex cones with $C \subseteq U$ and $D \subseteq V$ satisfying   $c(C,D) =c(U,V)$, $c(C,D) >c(U,V)$, and $c(C,D) <c(U,V)$, respectively. 
	In particular, these examples 
	illustrate that the similar inequality presented in \cref{lemma:CD}\cref{lemma:CD:UV} for minimal angle   does not hold for  angle. Suppose $\mathcal{H} =\mathbb{R}^{2}$.
	\begin{enumerate}
		\item \label{remark:c:=} Suppose $C:=\mathbb{R}_{+} (1,0) $, $D:=\mathbb{R}_{+} (-1,0) $, $U:=\mathbb{R}^{2}_{+}$  and $V:=\mathbb{R}^{2}_{-}$. Then $C \subseteq U$, $D \subseteq V$, $(C \cap D)^{\ominus}=\mathcal{H}$, $(U \cap V)^{\ominus} =\mathcal{H}$. Hence,  $c (C,D) =c_{0}(C,D)=0 =c_{0} (U,V)= c (U,V)$. 
		\item \label{remark:c:>} Suppose $C:=\mathbb{R}_{+} (1,0) $, $D:=\mathbb{R}_{+} (1,1) $, and $U=V=\mathcal{H}$.  Then $C \subseteq U$, $D \subseteq V$,  $(C \cap D)^{\ominus} =\mathcal{H}$, $(U \cap V)^{\ominus}=\{0\} $.  Hence,   $c (C,D) =c_{0}(C,D) =\frac{\sqrt{2}}{2} > 0 =c_{0}(\{0\}, \{0\}) = c (U,V)$.
		\item \label{remark:c:<} Suppose $C:=\mathbb{R}_{+} (1,0) $, $D:=\mathbb{R}_{+} (-1,0) $, 
$U:=\mathbb{R}^{2}_{+}$  and $V:=\{ (x_{1},x_{2}) ~:~  -x_{1} \geq x_{2} \geq 0 \}$.  Then $C \subseteq U$, $D \subseteq V$,  $U \cap V =\{0\}$. Hence,   $c (C,D) =c_{0}(C,D)=0 <\frac{\sqrt{2}}{2}= \Innp{(0, 1), (-\frac{\sqrt{2}}{2}, \frac{\sqrt{2}}{2})}=c_{0} (U,V)= c (U,V)$.
	\end{enumerate}
\end{remark}
	
\begin{fact} {\rm \cite[Propositions~3.3(4)]{DH2006II}} \label{fact:coneK1K2}
	Let  $K_{1}$ and $K_{2}$ be  nonempty  closed convex cones in $\mathcal{H}$. Then $c(K_{1},K_{2}) = c_{0} \left(  K_{1}  \cap (K_{1}\cap K_{2})^{\ominus}, K_{2} \cap (K_{1}\cap K_{2})^{\ominus} \right) $.
\end{fact}

\subsection*{Miscellaneous results}
The following fact is necessary for the proof of \cref{lemma:ineq}  below. 
\begin{fact}[Dunkl-Williams inequality]{\rm \cite[page~54]{DW1964}} \label{fact:DWineqH}
	Let $u$ and $v$ be in $\mathcal{H} \smallsetminus \{0\}$. Then 
	\begin{align*}
	\frac{1}{2} \Norm{ \frac{u}{\norm{u}} - \frac{v}{\norm{v}}  } \leq \frac{\norm{u-v}}{ \norm{u} +\norm{v}}.
	\end{align*}
\end{fact}
In fact, in the following result, $\gamma =\frac{\beta}{2}$ is motivated by \cite[Proposition~4.2]{SS2016I} which shows an identity on the maximal angle between two closed convex cones.
\begin{lemma} \label{lemma:ineq}
		 Let $C$ and $D$ be  cones in $\mathcal{H}$ such that $C \cap \mathbf{S}_{\mathcal{H}} \neq \varnothing$ and $D \cap \mathbf{S}_{\mathcal{H}} \neq \varnothing$.  Set
	\begin{align*}
	& \alpha:=  \sup \left\{   \innp{x,y} ~:~ x \in C \cap \mathbf{S}_{\mathcal{H}}, ~ y \in D \cap \mathbf{S}_{\mathcal{H}}\right\},\\
	& \beta:=\dist \left( C \cap \mathbf{S}_{\mathcal{H}}, D \cap \mathbf{S}_{\mathcal{H}} \right) = \inf \left\{   \norm{x-y} ~:~ x \in C \cap \mathbf{S}_{\mathcal{H}}, ~ y \in D \cap \mathbf{S}_{\mathcal{H}}\right\},\\
	&\gamma:= \inf \left\{   \frac{ \norm{x-y}} {\norm{x} +\norm{y}}~:~ x \in C, ~ y \in D, ~(x,y) \neq (0,0) \right\}.
	\end{align*}	
	Then $\alpha  =1 -\frac{\beta^{2}}{2} =1-2 \gamma^{2}$.
\end{lemma}

\begin{proof}
 It is easy to see that
\begin{align*}
\frac{ \beta^{2} }{2} & = \frac{ 1}{2}  \inf \left\{   \norm{x-y}^{2} ~:~ x \in C \cap \mathbf{S}_{\mathcal{H}}, ~ y \in D \cap \mathbf{S}_{\mathcal{H}}\right\}\\
& =\frac{ 1}{2}   \inf \left\{  2 -2\innp{x,y} ~:~ x \in C \cap \mathbf{S}_{\mathcal{H}}, ~ y \in D \cap \mathbf{S}_{\mathcal{H}}\right\}\\
& =1+  \inf \left\{  - \innp{x,y} ~:~ x \in C \cap \mathbf{S}_{\mathcal{H}}, ~ y \in D \cap \mathbf{S}_{\mathcal{H}}\right\} \\
& =1 -  \sup \left\{  \innp{x,y} ~:~ x \in C \cap \mathbf{S}_{\mathcal{H}}, ~ y \in D \cap \mathbf{S}_{\mathcal{H}}\right\} =1-\alpha,
\end{align*}
which implies that $ \alpha =1 -\frac{\beta^{2}}{2} $.

It remains to show that $\gamma =  \frac{\beta}{2}$. Clearly,
\begin{align*}
\gamma & =	\inf \left\{   \frac{ \norm{x-y}} {\norm{x} +\norm{y}}~:~ x \in C, ~ y \in D,~ (x,y) \neq (0,0) \right\} \\
&\leq 	\inf \left\{   \frac{ \norm{x-y}} {\norm{x} +\norm{y}}~:~ x \in C \cap  \mathbf{S}_{\mathcal{H}},  ~ y \in D \cap  \mathbf{S}_{\mathcal{H}} \right\}\\
&	= \inf \left\{   \frac{ \norm{x-y}} {2}~:~ x \in C \cap  \mathbf{S}_{\mathcal{H}},  ~ y \in D \cap  \mathbf{S}_{\mathcal{H}} \right\}=\frac{\beta}{2}.
\end{align*}

On the other hand, let $x \in C$ and  $ y \in D$ such that $(x,y) \neq (0,0)$. Because for every $z \in C \cap \mathbf{S}_{\mathcal{H}}$ and   $w \in D \cap \mathbf{S}_{\mathcal{H}}$, $\norm{z-w} \leq \norm{z}+\norm{w}=2$, we have that $\beta \leq 2$.
If $x=0$ and $y \neq 0$, then $ \frac{ \norm{x-y}} {\norm{x} +\norm{y}} =1 \geq \frac{\beta}{2}$. Similarly, if $x\neq 0$ and $y = 0$, then $ \frac{ \norm{x-y}} {\norm{x} +\norm{y}} =1 \geq \frac{\beta}{2}$.
Assume $x\neq 0$ and $y \neq 0$. Substitute  $u=x$ and $v=y$ in \cref{fact:DWineqH} to see that
\begin{align*}
\frac{\norm{x-y}}{ \norm{x} +\norm{y}} \geq 	\frac{1}{2} \Norm{ \frac{x}{\norm{x}} - \frac{y}{\norm{y}}  } \geq 	\frac{1}{2} \inf \left\{   \norm{a-b} ~:~ a \in C  \cap \mathbf{S}_{\mathcal{H}}, ~ b \in D \cap \mathbf{S}_{\mathcal{H}}\right\}= \frac{\beta}{2}.
\end{align*}
Combine these results to see that  $\gamma \geq  \frac{\beta}{2}$. Hence, we obtain that $\gamma =  \frac{\beta}{2}$.
\end{proof}

\section{Principal vectors of  convex sets} \label{sec:Principalvectors}

 In this section, we shall provide equivalent expressions of the cosine of the minimal angle of two convex sets. Moreover, we shall also construct necessary conditions for one pair of vectors in $\mathcal{H}\times \mathcal{H}$ to be a pair of principal vectors defined in \cref{defn:principalvector}, which is actually one optimal solution for the related minimal angle problem.

\begin{theorem} \label{theorem:c0}
 Let $C$ and $D$ be nonempty   convex subsets of $\mathcal{H}$  such that    $C \neq \{0\}$ and $D \neq \{0\}$.  Set
	\begin{align*}
	& \alpha:=  \sup \left\{   \innp{x,y} ~:~ x \in \overline{\cone}(C) \cap \mathbf{S}_{\mathcal{H}}, ~ y \in \overline{\cone}(D) \cap \mathbf{S}_{\mathcal{H}}\right\},\\
	& \beta:=\dist \left( \overline{\cone}(C)  \cap \mathbf{S}_{\mathcal{H}}, \overline{\cone}(D)  \cap \mathbf{S}_{\mathcal{H}} \right) = \inf \left\{   \norm{x-y} ~:~ x \in \overline{\cone}(C) \cap \mathbf{S}_{\mathcal{H}}, ~ y \in \overline{\cone}(D) \cap \mathbf{S}_{\mathcal{H}}\right\},\\
	&\gamma:= \inf \left\{   \frac{ \norm{x-y}} {\norm{x} +\norm{y}}~:~ x \in \overline{\cone}(C), ~ y \in \overline{\cone}(D), ~(x,y) \neq (0,0) \right\}.
	\end{align*}	
	Then the following statements hold:
	\begin{enumerate}
		\item \label{theorem:c0:P} $c_{0}(C,D)   =\sup \left\{  \sqrt{ \innp{\Pro_{\overline{\cone}(C)}x,\Pro_{\overline{\cone}(D)}\Pro_{\overline{\cone}(C)}x} }~:~  x \in  \mathbf{S}_{\mathcal{H}} \right\}$.
		\item \label{theorem:c0:neq0} Assume that $c_{0}(C,D) \neq 0$. 
		Then $c_{0}(C,D) =\alpha  =1 -\frac{\beta^{2}}{2} =1-2 \gamma^{2}$.
	\end{enumerate}
\end{theorem}

\begin{proof}
	\cref{theorem:c0:P}:  Note that by \cref{fact:chara:PK}, 
	\begin{align}\label{eq:theorem:c0:P}
	\innp{ \Pro_{\overline{\cone}(C)} x  - \Pro_{\overline{\cone}(D)}\Pro_{\overline{\cone}(C)} x,  \Pro_{\overline{\cone}(D)}\Pro_{\overline{\cone}(C)} x  }=0.
	\end{align} 
	Moreover, by \cite[Proposition~2.2(2), Theorem~2.5]{DH2006II} and \cite[Theorem~5.6(7)]{D2012}, we know that
	\begin{align*}
c_{0}(C,D) &~=~ \sup \left\{ \norm{ \Pro_{\overline{\cone}(D)}\Pro_{\overline{\cone}(C)}x }~:~  x \in  \mathbf{S}_{\mathcal{H}} \right\}\\
&~=~\sup \left\{  \sqrt{ \innp{ \Pro_{\overline{\cone}(D)}\Pro_{\overline{\cone}(C)}x ,\Pro_{\overline{\cone}(D)}\Pro_{\overline{\cone}(C)}x }  }~:~  x \in  \mathbf{S}_{\mathcal{H}} \right\}\\
&\stackrel{\cref{eq:theorem:c0:P}}{=} \sup  \left \{  \sqrt{ \innp{ \Pro_{\overline{\cone}(C)} x,\Pro_{\overline{\cone}(D)}\Pro_{\overline{\cone}(C)} x}  }~:~  x \in  \mathbf{S}_{\mathcal{H}} \right \}.
	\end{align*}

	\cref{theorem:c0:neq0}: 	
	By \cref{fact:coBS},  $	c_{0}(C,D)= \sup \left\{   \innp{x,y} ~:~ x \in \overline{\cone} (C) \cap \mathbf{S}_{\mathcal{H}}, ~ y \in \overline{\cone} (D) \cap \mathbf{S}_{\mathcal{H}}\right\}$, since  $c_{0}(C,D) \neq 0$.
	Moreover, $C$ and $D$ being nonempty, $C \neq \{0\}$ and $D \neq \{0\}$ imply that $\overline{\cone}(C) \cap \mathbf{S}_{\mathcal{H}} \neq \varnothing$ and  $  \overline{\cone}(D) \cap \mathbf{S}_{\mathcal{H}} \neq \varnothing$.
	Hence, apply \cref{lemma:ineq} with $C$ and $D$ replaced by $ \overline{\cone} (C)$ and $  \overline{\cone} (D)$ respectively  to see that  \cref{theorem:c0:neq0} is true.
\end{proof}

\begin{definition} \label{defn:principalvector} 
	Let $C$ and $D$ be  nonempty convex subsets of $\mathcal{H} $.  Let $(\bar{x}, \bar{y}) \in \mathcal{H} \times \mathcal{H}$. We say that $(\bar{x}, \bar{y})$ is a pair of  \emph{principal vectors} of $C$ and $D$, if  $\bar{x} \in \overline{\cone}(C) \cap \mathbf{B}_{\mathcal{H}}$,  $\bar{y} \in \overline{\cone}(D)  \cap \mathbf{B}_{\mathcal{H}}$, and 
	\begin{align*}
	\innp{ \bar{x}, \bar{y}} 
=	\sup \left\{   \innp{x,y} ~:~ x \in \overline{\cone}(C)  \cap \mathbf{B}_{\mathcal{H}}, ~ y \in \overline{\cone}(D)  \cap \mathbf{B}_{\mathcal{H}} \right\}.
	\end{align*}
\end{definition}

Note that a pair of principal vectors of convex sets $C$ and $D$ is an optimal solution of the minimal angle problem associated with $C$ and $D$, and    that when $C$ and $D$ are nonzero linear subspaces, 
 it is easy to see that $(\bar{x}, \bar{y})$ is a pair of  principal vectors of $C$ and $D$, if  $\bar{x} \in C \cap \mathbf{S}_{\mathcal{H}}$,  $\bar{y} \in D  \cap \mathbf{S}_{\mathcal{H}}$, and $\innp{ \bar{x}, \bar{y}} = 	\sup \left\{   \innp{x,y} ~:~ x \in C  \cap \mathbf{S}_{\mathcal{H}}, ~ y \in D  \cap \mathbf{S}_{\mathcal{H}} \right\}$.

The following result provides  sufficient conditions for the existence of optimal solutions of the minimal angle problem associated with two nonempty and nonzero convex subsets of the Hilbert space.  In particular,  these optimal solutions   always exist in finite-dimensional spaces. 
\begin{proposition} \label{prop:optimalpair:exist}
	 Let $\mathcal{K}$ be a finite-dimensional linear subspace of $\mathcal{H}$.  	Let $C$ and $D$ be nonempty convex subsets of $\mathcal{H}$ such that  $C \neq \{0\}$ and $D \neq \{0\}$.  Suppose that $C \subseteq  \mathcal{K}$ or $D \subseteq  \mathcal{K}$.  
	Then there exists $(\bar{x}, \bar{y}) \in \mathcal{H} \times \mathcal{H}$ such that $(\bar{x}, \bar{y})$ is a pair of principal vectors of $C$ and $D$. 
\end{proposition}

\begin{proof}
	By \cref{fact:AnglesProperties}\cref{fact:AnglesProperties:between0and1},  
	we have exactly the following two cases.
	
	\emph{Case~1:} $c_{0} (C, D ) =0$.  Then it is clear that $(\bar{x}, \bar{y})=(0,0)$ works for this special case.
	
	\emph{Case~2:}  $c_{0} (C, D ) >0$. Then by \cref{fact:coBS}, $	c_{0} (C, D )=\sup \left\{   \innp{x,y} ~:~ x \in \overline{\cone}(C)  \cap \mathbf{S}_{\mathcal{H}}, ~ y \in \overline{\cone}(D)  \cap \mathbf{S}_{\mathcal{H}} \right\}$.
	Hence,	there exist sequences $(x_{k})_{k \in \mathbb{N}}$ in $\overline{\cone}(C)   \cap \mathbf{S}_{\mathcal{H}}$ and $ (y_{k})_{k \in \mathbb{N}} $ in $\overline{\cone}(D)  \cap \mathbf{S}_{\mathcal{H}}$ such that $ \innp{x_{k} ,y_{k}} \to c_{0} (C, D )$.
	Now, because $(x_{k})_{k \in \mathbb{N}}$ and $ (y_{k})_{k \in \mathbb{N}} $ are in $\mathbf{S}_{\mathcal{H}}$, and $\overline{\cone}(C) $ and $\overline{\cone}(D)$ are nonempty closed and convex,
	by \cite[Lemma~2.45 and Corollary~3.35]{BC2017} and by passing to a subsequence if necessary, there exist $\bar{x} \in \overline{\cone}(C)$ and $\bar{y} \in \overline{\cone}(D)$ such that $x_{k} \weakly \bar{x}$ and $y_{k} \weakly \bar{y}$. By \cite[Lemma~2.42]{BC2017}, we have that 
	\begin{align*} 
	\norm{\bar{x} } \leq \liminf_{k \to \infty} \norm{x_{k}}  =1 \quad \text{ and} \quad  \norm{\bar{y} }   \leq \liminf_{k \to \infty} \norm{y_{k}} =1,
	\end{align*}
	which implies that $\bar{x} \in \overline{\cone}(C) \cap \mathbf{B}_{\mathcal{H}}$ and  $\bar{y} \in \overline{\cone}(D) \cap \mathbf{B}_{\mathcal{H}}$.
	
	Moreover, by assumption, without loss of generality, $\overline{\cone}(C) \subseteq \mathcal{K}$. Then  $x_{k} \weakly \bar{x}$ in $\mathcal{K}$ and so, by \cite[Lemma~2.51(ii)]{BC2017}, $x_{k} \to \bar{x}$. Hence, by \cite[Lemma~2.51(iii)]{BC2017}, $x_{k} \to \bar{x}$ and $y_{k} \weakly \bar{y}$ imply that $  \innp{x_{k}, y_{k}} \to 	\innp{\bar{x}, \bar{y}} $.   Combine this with the result $ \innp{x_{k} ,y_{k}} \to c_{0} (C, D )$ to deduce the required result.
\end{proof}

\begin{lemma} \label{lemma:optimalpair:norm1}
	Let $C$ and $D$ be nonempty  convex subsets of $\mathcal{H}$ with $C \neq \{0\}$ and $D \neq \{0\}$.  Assume  that $(\bar{x}, \bar{y})$ is a pair of principal vectors  of  $C$ and $D$ and that $\innp{ \bar{x}, \bar{y}}  \neq 0$. Then $\norm{\bar{x} } = 1$ and $\norm{\bar{y} } = 1$.   		
	Consequently, $\innp{ \bar{x}, \bar{y}} = \sup \left\{   \innp{x,y} ~:~ x \in \overline{\cone}(C) \cap \mathbf{S}_{\mathcal{H}}, ~ y \in \overline{\cone}(D) \cap \mathbf{S}_{\mathcal{H}}\right\} $.	
\end{lemma}

\begin{proof}
	Because  $\innp{ \bar{x}, \bar{y}}  \neq 0$ and $0 \in \overline{\cone}(C) \cap \overline{\cone}(D) \cap \mathbf{B}_{\mathcal{H}}$, we know that 
	\begin{align}\label{eq:lemma:optimalpair:norm1}
	\innp{ \bar{x}, \bar{y}}  > 0,
	\end{align}
	which implies that $\bar{x} \neq 0$ and $ \bar{y} \neq 0$. Suppose to the contrary that $\norm{\bar{x} } \neq  1$ or $\norm{\bar{y} } \neq 1$.    Then by $\bar{x} \in \overline{\cone}(C) \cap \mathbf{B}_{\mathcal{H}}$ and $\bar{y} \in \overline{\cone}(D) \cap \mathbf{B}_{\mathcal{H}}$, we have that $1<\frac{1}{ \norm{\bar{x}}\norm{\bar{y} }} $ and that 
	\begin{align}\label{eq:lemma:optimalpair:norm1:in}
	\frac{\bar{x} }{\norm{\bar{x}}} \in  \overline{\cone}(C) \cap \mathbf{B}_{\mathcal{H}}\quad  \text{and} \quad \frac{\bar{y} }{\norm{\bar{y}}} \in  \overline{\cone}(D) \cap \mathbf{B}_{\mathcal{H}}.
	\end{align}
In addition,
	\begin{align*}
	\Innp{ 	\frac{\bar{x}}{ \norm{\bar{x}}}, \frac{\bar{y}}{ \norm{\bar{y}} } } =	\frac{1}{ \norm{\bar{x}}\norm{\bar{y} }} \innp{ \bar{x}, \bar{y}} \stackrel{\cref{eq:lemma:optimalpair:norm1}}{>}  \innp{ \bar{x}, \bar{y}}=\sup \left\{   \innp{x,y} ~:~ x \in \overline{\cone}(C) \cap \mathbf{B}_{\mathcal{H}}, ~ y \in \overline{\cone}(D) \cap \mathbf{B}_{\mathcal{H}} \right\},
	\end{align*}
	which contradicts  \cref{eq:lemma:optimalpair:norm1:in}. Therefore, $\norm{\bar{x} } = 1$ and $\norm{\bar{y} } = 1$.  
\end{proof}
The idea of the following proof in  case 2  is from that of \cite[Proposition~1.3]{SS2016I} which shows   necessary conditions for  optimal solutions of the maximization problem defined in  \cite[Definition~1.1]{SS2016I}. 

\begin{lemma} \label{lemma:K1K2polar}
	Let $C$ and $D$ be nonempty  convex subsets of $\mathcal{H}$ with $C \neq \{0\}$ and $D \neq \{0\}$. Assume  that $(\bar{x}, \bar{y})$ is a  pair of principal vectors  of  $C$ and $D$.  Then
	\begin{align*}
	\bar{y} - \innp{\bar{x},\bar{y} } \bar{x} \in C^{\ominus}  \quad \text{and}  \quad \bar{x} - \innp{\bar{x},\bar{y} } \bar{y} \in D^{\ominus}. 
	\end{align*}
\end{lemma}
\begin{proof}	
	By \cref{defn:principalvector} and \cref{fact:AnglesProperties}\cref{fact:AnglesProperties:between0and1},   $\innp{\bar{x}, \bar{y}}=c_{0} (C,D) \in \left[ 0,1\right]$. Hence,
	we have exactly the following two cases:
	
	\emph{Case~1:} $c_{0} (C,D)=0$.   Then by \cref{fact:AnglesProperties}\cref{fact:AnglesProperties:c0:EQ}, $c_{0} \left( \overline{\cone} (C) ,\overline{\cone} (D) \right)=c_{0} (C,D)=0$. Hence, by 	 \cite[Theorem~5.8]{DH2006II}  and \cref{fact:dualcone}\cref{fact:dualcone:eq},
	\begin{align*}
	&\bar{y} - \innp{\bar{x},\bar{y} } \bar{x} = 	\bar{y}  \in \overline{\cone} (D)  \subseteq (\overline{\cone} (C) )^{\ominus} = C ^{\ominus}, \\
&\bar{x} - \innp{\bar{x},\bar{y} } \bar{y} = \bar{x}  \in \overline{\cone} (C) \subseteq (\overline{\cone} (D))^{\ominus} =D^{\ominus}. 
	\end{align*}
	
	\emph{Case~2:}  $c_{0} (C,D)>0$.  Then by \cref{lemma:optimalpair:norm1}, $\bar{x}  \in \overline{\cone} (C)  \cap \mathbf{S}_{\mathcal{H}}$,  $\bar{y}  \in \overline{\cone} (D)  \cap \mathbf{S}_{\mathcal{H}}$, and 
	\begin{align*}
	\innp{\bar{x},\bar{y} }   =\sup \left\{   \innp{x,y} ~:~ x \in \overline{\cone} (C)  \cap \mathbf{S}_{\mathcal{H}}, ~ y \in \overline{\cone} (D)  \cap \mathbf{S}_{\mathcal{H}} \right\},
	\end{align*}
	which implies that 
	\begin{align}\label{eq:lemma:K1K2polar}
	( \forall x \in \overline{\cone} (C)  \cap \mathbf{S}_{\mathcal{H}} ) \quad \innp{\bar{x},\bar{y} } \geq \innp{x,\bar{y} }.
	\end{align}
	Let $x \in \overline{\cone} (C)  \smallsetminus \{0\}$. Set $\epsilon_{x}:=\frac{\norm{\bar{x}}}{\norm{x}}$. Note that 
	\begin{align*}
	(\forall t \in \left[0, \epsilon_{x}\right[\,) \quad 	\norm{\bar{x} +tx}\geq \norm{\bar{x} } -t\norm{x} >0. 
	\end{align*}
	Hence, the function $f: \left[0, \epsilon_{x}\right[ \to \mathbb{R} : t \mapsto \Innp{ \frac{ \bar{x} +tx}{	\norm{\bar{x} +tx}}, \bar{y}}$ is well-defined. Because $(\forall t \in  \left[0, \epsilon_{x}\right[\, )$ $f(t) = \frac{ \innp{\bar{x}, \bar{y}} +t\innp{x,\bar{y}}}{	\norm{\bar{x} +tx}}$, by \cite[Example~2.65]{BC2017}, we have that
	\begin{align} \label{eq:lemma:K1K2polar:deri}
	(\forall t \in  \left[0, \epsilon_{x}\right[\, ) \quad f^{\prime}_{+} (t) =  \frac{ 1 }{	\norm{\bar{x} +tx}^{2}} \left(  \innp{x, \bar{y}}\norm{\bar{x} +tx} - (   \innp{\bar{x}, \bar{y}} +t\innp{x,\bar{y}} )  \Innp{\frac{ \bar{x} +tx}{	\norm{\bar{x} +tx}},x} \right).
	\end{align}
	Because $\{ \bar{x}, x \} \subseteq \overline{\cone} (C)$ and $\overline{\cone} (C)$ is a closed convex cone, by \cref{lemma:K}, $(\forall t \in  \left[0, \epsilon_{x}\right[\,)$ $ \frac{ \bar{x} +tx}{	\norm{\bar{x} +tx}} \in \overline{\cone} (C)+\overline{\cone} (C) = \overline{\cone} (C)$. Hence,   by \cref{eq:lemma:K1K2polar}, $f(0) = \max \{ f(t) ~:~t \in \left[0, \epsilon_{x}\right[ \}$. Combine this with \cref{eq:lemma:K1K2polar:deri} to see that
	\begin{align*}
	0 \geq  f^{\prime}_{+} (0) =  \innp{x, \bar{y}} - \innp{\bar{x}, \bar{y}} \innp{\bar{x}, x} = \Innp{\bar{y} - \innp{\bar{x}, \bar{y}} \bar{x}, x }.
	\end{align*}
	Therefore, $(\forall x \in \overline{\cone} (C))$ $  \Innp{\bar{y} - \innp{\bar{x}, \bar{y}} \bar{x},x } \leq 0$, that is, $	\bar{y} - \innp{\bar{x},\bar{y} } \bar{x} \in( \overline{\cone} (C))^{\ominus} =C^{\ominus}$.
	
	By similar arguments, we get  $\bar{x} - \innp{\bar{x},\bar{y} } \bar{y} \in (\overline{\cone} (D))^{\ominus}=D^{\ominus}$.
\end{proof}

\cite[Theorem~2.3]{SS2016I} considers points on the boundary of closed convex cones relative to a linear subspace. 
The idea of the following proof is from \cite[Theorem~2.3]{SS2016I}.
\begin{lemma}\label{lemma:K1K2incomple}
	Let $C$ and $D$ be nonempty subsets of $\mathcal{H}$.  Assume  that $x \in \mathbf{S}_{\mathcal{H}}$, $y  \in \mathbf{S}_{\mathcal{H}}$, $y -\innp{x,y}x \in C^{\ominus}$, and $x-\innp{x,y}y \in D^{\ominus}$, and that    $\abs{ \innp{x,y} }\neq 1$. Then 
	\begin{align*}
x +\cone (y - \innp{x,y } x ) \subseteq  C^{c},  ~ x \notin \inte C,  ~
y+\cone ( x - \innp{x,y} y ) \subseteq D^{c},  \text{ and } y \notin \inte D.
	\end{align*}
	
\end{lemma}

\begin{proof}
	Note that  by \cref{fact:coneC}, $  x +\cone (y - \innp{x,y } x ) = x + \mathbb{R}_{++} (y - \innp{x,y } x ) $.
	
	Let $\alpha \in \mathbb{R}_{++}$.  Assume to the contrary that $x+\alpha (y - \innp{x,y} x ) \in C$. By  the assumption $\abs{ \innp{x,y} } \neq 1$, we know that 
	\begin{align}\label{eq:coro:K1K2incomple:>0}
	\norm{y - \innp{x,y} x}^{2} =1-\innp{x,y}^{2} >0.
	\end{align}
	Note that $ \innp{ x, 	y - \innp{x,y} x  } =\innp{x, 	y } -\innp{x, 	y} \norm{x}^{2}=0$.
	Combine the assumptions that $y -\innp{x,y}x \in C^{\ominus}$ and $x+\alpha (y - \innp{x,y} x ) \in C$ with    \cref{eq:coro:K1K2incomple:>0} to  obtain that
	\begin{align*}
	0 \geq \Innp{x+\alpha (y - \innp{x,y} x ) , y -\innp{x,y}x   } = \alpha \Innp{y - \innp{x,y} x, y -\innp{x,y}x } = \alpha \norm{y -\innp{x,y}x}^{2} >0,
	\end{align*}	
	which is a contradiction. So, $x+\alpha (y - \innp{x,y} x ) \in C^{c} $. Because $\alpha \in \mathbb{R}_{++}$ is arbitrary, we know that  $  x +\cone (y - \innp{x,y } x ) = x + \mathbb{R}_{++} (y - \innp{x,y } x ) \subseteq  C^{c}$, and that  $x  \notin \inte C$. 
	
	By analogous arguments, we get $y+\cone ( x - \innp{x,y} y ) \subseteq D^{c}$ and  $y \notin \inte D$ as well. 
\end{proof}

\begin{lemma} \label{lemma:xy}
	Let $C$  be a nonempty convex subset of $\mathcal{H}$.  Let $(x, y) \in \mathcal{H} \times  \mathcal{H}$.   Assume that $x \in \mathbf{S}_{\mathcal{H}}$.  Then the following equivalences are true:
	\begin{enumerate}
		\item \label{lemma:xy:o}  $ y -\innp{x,y} x \in C^{\ominus} $ $\Leftrightarrow$ $y -\innp{x,y} x \in (C - x)^{\ominus} $.
		\item \label{lemma:xy:EQ} $x = \Pro_{C} \left(x + \cone(y -\innp{x,y} x)  \right)$ $\Leftrightarrow$  $ [ $ $x \in C$ and $ y -\innp{x,y} x \in C^{\ominus} $ $]$.
		\item \label{lemma:xy:cone:EQ}   $(\forall \lambda \in \mathbb{R}_{++} )$ $\lambda( y -\innp{x,y} x)= \Pro_{C^{\ominus}} (x + \lambda (y -\innp{x,y} x )) $ $\Leftrightarrow$  $ [ $ $x \in \overline{\cone}(C)$ and $ y -\innp{x,y} x \in C^{\ominus} $ $]$. 
	\end{enumerate} 
\end{lemma}
\begin{proof}
	Because $x \in \mathbf{S}_{\mathcal{H}}$,  we have  
	\begin{align}\label{eq:lemma:xy}
	\innp{ y -\innp{x,y} x , x }=\innp{y,x} -\innp{x,y}\norm{x}^{2}=0.
	\end{align}
	
	\cref{lemma:xy:o}: 	By the definition of polar cone, 	
	\begin{align*}
	y -\innp{x,y} x \in C^{\ominus}  & ~\Leftrightarrow~ (\forall z \in C) ~ \innp{y -\innp{x,y} x,z} \leq 0\\
	&\stackrel{\cref{eq:lemma:xy}}{\Leftrightarrow}   (\forall z \in C)~  \innp{y -\innp{x,y} x, z-x} \leq 0\\
	& ~\Leftrightarrow~ 	y -\innp{x,y} x \in (C -x)^{\ominus}. 
	\end{align*}
	
	\cref{lemma:xy:EQ}: Because, by  \cref{fact:coneC}, $\cone(y -\innp{x,y} x)=\mathbb{R}_{++} (y -\innp{x,y} x)$, we have that 
	\begin{align*}
	&~x = \Pro_{C} \left(x + \cone(y -\innp{x,y} x)  \right) \\
	\Leftrightarrow & ~(\forall \lambda \in \mathbb{R}_{++} )~	x = \Pro_{C} \left(x + \lambda(y -\innp{x,y} x)  \right)\\
	\Leftrightarrow 	 &~\left[x \in C \text{ and } (\forall \lambda \in \mathbb{R}_{++} ) (\forall z \in C) ~\Innp{ x + \lambda(y -\innp{x,y} x)  -x, z-x } \leq 0 \right] \quad (\text{by \cref{fact:characterPC}})\\
	\Leftrightarrow & ~ [x \in C \text{ and } (\forall z \in C)~ \innp{  y -\innp{x,y} x, z-x } \leq 0 ] \quad (\text{by $\lambda \in \mathbb{R}_{++}$})\\
	\Leftrightarrow & ~ [x \in C \text{ and } y -\innp{x,y} x \in (C -x)^{\ominus}]\\
	\stackrel{\text{\cref{lemma:xy:o}}}{\Leftrightarrow} &~ [x \in C \text{ and } y -\innp{x,y} x \in C^{\ominus}].
	\end{align*}
	
	\cref{lemma:xy:cone:EQ}:  Let $\lambda \in \mathbb{R}_{++}$. 
	Using  \cref{fact:characterPC} and \cref{fact:dualcone}\cref{fact:dualcone:o}, we have that 
	\begin{align*}
	&\lambda ( y -\innp{x,y} x)= \Pro_{C^{\ominus}} \left(x +\lambda (y -\innp{x,y} x) \right) \\
	~\Leftrightarrow~ &y -\innp{x,y} x \in C^{\ominus}  \text{ and }   (\forall z \in C^{\ominus} )~\Innp{x +  \lambda ( y -\innp{x,y} x) - \lambda ( y -\innp{x,y} x), z- \lambda(y -\innp{x,y} x)} \leq 0\\
~\Leftrightarrow~&y -\innp{x,y} x \in C^{\ominus}  \text{ and } (\forall z \in C^{\ominus} )~\Innp{x, z-\lambda (y -\innp{x,y} x)} \leq 0\\
	\stackrel{\cref{eq:lemma:xy}}{\Leftrightarrow} ~& y -\innp{x,y} x \in C^{\ominus}  \text{ and } x \in (C^{\ominus} )^{\ominus}\\
~\Leftrightarrow~&y -\innp{x,y} x \in C^{\ominus}  \text{ and } x \in \overline{\cone}C,
	\end{align*}
	where the last equivalence is from  \cref{fact:dualcone}\cref{fact:dualcone:oo}.
\end{proof}
The following result provides necessary conditions for  $(x ,y) \in \mathcal{H} \times \mathcal{H}$   to be a pair of  principal vectors  of two nonempty  convex  subsets  $C$ and $D$ in $\mathcal{H}$ with $C \neq \{0\}$ and $D \neq \{0\}$. 

\begin{proposition} \label{prop:optimalpair:neces}
	Let $C$ and $D$ be nonempty  convex subsets of $\mathcal{H}$ with $C \neq \{0\}$ and $D \neq \{0\}$.  Assume  that $(x ,y)$ is a pair of principal vectors  of  $C$ and $D$. Let $ \lambda \in \mathbb{R}_{++}$. Then the following statements hold:
	\begin{enumerate}
		\item \label{prop:optimalpair:neces:inte} If $x\neq 0$, $y\neq 0$ and  $\abs{\innp{x,y}} \neq 1$, then $x +\cone (y - \innp{x,y } x ) \subseteq  (\overline{\cone}(C))^{c}$,   $x \notin \inte \overline{\cone}(C)$, $y+\cone ( x - \innp{x,y} y ) \subseteq (\overline{\cone}(D))^{c}$, $y \notin \inte \overline{\cone}(D)$.
 
		\item \label{prop:optimalpair:neces:proje} $x = \Pro_{\overline{\cone}(C)} \left(x + \cone(y -\innp{x,y} x)  \right) $,   and  $y = \Pro_{\overline{\cone}(D)} \left( y +  \cone(x -\innp{x,y}y)\right) $.
		\item \label{prop:optimalpair:neces:o} $\lambda (y -\innp{x,y}x)= \Pro_{C^{\ominus}} \left( x+ \lambda( y -\innp{x,y}x)\right) $,   and $\lambda (x -\innp{x,y}y)= \Pro_{D^{\ominus}}  ( y +\lambda (x -\innp{x,y}y) ) $.
	\end{enumerate}
\end{proposition}
\begin{proof}
	By \cref{defn:principalvector} and \cref{fact:AnglesProperties}\cref{fact:AnglesProperties:between0and1},   $\innp{x,y}=c_{0} (C,D) \in \left[ 0,1\right]$. Hence, we have exactly the following two cases:
	
	\emph{Case~1:} $\innp{x,y}=0$. 
	By \cref{defn:principalvector},  \cref{lemma:K1K2polar} and \cref{fact:dualcone}\cref{fact:dualcone:eq}, $x \in \overline{\cone}(C) \cap D^{\ominus}=\overline{\cone}(C) \cap (\overline{\cone}(D))^{\ominus}$ and $y \in \overline{\cone}(D) \cap C^{\ominus}= \overline{\cone}(D) \cap (\overline{\cone}(C))^{\ominus}$.

	For the proof of \cref{prop:optimalpair:neces:inte}, because $\innp{x,y}=0$ and $y\neq 0$, we have that 
	\begin{align*}
	(\forall \alpha \in \mathbb{R}_{++}) \quad    \innp{ x +\alpha (y - \innp{x,y } x ), y} =\innp{x,y}+\alpha\innp{y,y}=\alpha\norm{y}^{2} >0,  
	\end{align*}
	which, combining with $y \in (\overline{\cone}(C))^{\ominus} $, implies that $x +\cone (y - \innp{x,y } x ) =x +\mathbb{R}_{++} (y - \innp{x,y } x ) \subseteq  (\overline{\cone}(C))^{c}$ and that $x \notin \inte \overline{\cone}(C)$.
	
	Consider the proof of \cref{prop:optimalpair:neces:proje}$\&$\cref{prop:optimalpair:neces:o}. Because $y \in (\overline{\cone}(C))^{\ominus}=C^{\ominus}$ and $x \in  \overline{\cone}(C) $,   we have that 
	\begin{subequations}
		\begin{align}
		&(\forall z \in  \overline{\cone}(C)) \quad \innp{ x+\lambda y -x, z-x} =\lambda \innp{y,z-x} =\lambda \innp{y,z} \leq 0, \label{eq:prop:optimalpair:neces:leq0}\\
		&(\forall w \in C^{\ominus}) \quad \innp{ x+\lambda y -\lambda y , w-\lambda y } =\innp{ x , w-\lambda y } = \innp{x,w} \leq 0.\label{eq:prop:optimalpair:neces:o}
		\end{align}
	\end{subequations}
Moreover, by \cref{eq:prop:optimalpair:neces:leq0}, $\innp{x,y} =0$  and  \cref{fact:characterPC}, we see that $x = \Pro_{ \overline{\cone}(C) } \left(x + \cone(y -\innp{x,y} x)  \right) $. In addition, by \cref{eq:prop:optimalpair:neces:o}  and  \cref{fact:characterPC}, we obtain that 
$\lambda (y -\innp{x,y}x)= \Pro_{C^{\ominus}} \left( x+ \lambda( y -\innp{x,y}x)\right) $.
	
	By similar arguments, we get the remaining parts of \cref{prop:optimalpair:neces:inte}, 	\cref{prop:optimalpair:neces:proje} and \cref{prop:optimalpair:neces:o}. 	
	
	\emph{Case~2:} $\innp{x,y} > 0$. Then by \Cref{lemma:optimalpair:norm1,lemma:K1K2polar} and \cref{fact:dualcone}\cref{fact:dualcone:eq},  we have that $ x \in \overline{\cone} (C) \cap \mathbf{S}_{\mathcal{H}}$, $ y \in \overline{\cone} (D)  \cap \mathbf{S}_{\mathcal{H}}$,   $y -\innp{x,y}x \in C^{\ominus}  =(\overline{\cone} (C) )^{\ominus} $, and $x -\innp{x,y}y  \in D^{\ominus}=(\overline{\cone} (D) )^{\ominus} $.  
	
	Now, apply \cref{lemma:K1K2incomple} with $C=\overline{\cone} (C) $ and $D=\overline{\cone} (D) $  to obtain the desired results in \cref{prop:optimalpair:neces:inte}. 
	
	Because $x \in \overline{\cone} (C) \cap \mathbf{S}_{\mathcal{H}}$ and    $ y -\innp{x,y}x \in (\overline{\cone} (C) )^{\ominus}  $, apply  \cref{lemma:xy}\cref{lemma:xy:EQ}$\&$\cref{lemma:xy:cone:EQ} with $C=\overline{\cone} (C) $ to see that $x = \Pro_{\overline{\cone} (C) } \left(x + \cone(y -\innp{x,y} x)  \right)$ and  $\lambda (y -\innp{x,y}x)= \Pro_{C^{\ominus}} \left( x+ \lambda( y -\innp{x,y}x)\right)$. 
	
	Similarly,  using $ y \in \overline{\cone} (D)  \cap \mathbf{S}_{\mathcal{H}}$ and  $x -\innp{x,y}y  \in (\overline{\cone} (D) )^{\ominus}  $, and applying  \cref{lemma:xy}\cref{lemma:xy:EQ}$\&$\cref{lemma:xy:cone:EQ}  with $C=\overline{\cone} (D) $ and with switching $x$ and $y$, we obtain that 
	$y = \Pro_{\overline{\cone} (D) } \left( y +  \cone(x -\innp{x,y}y)\right)  $ and $\lambda (x -\innp{x,y}y)= \Pro_{D^{\ominus}}   (y +\lambda (x -\innp{x,y}y) ) $.
	
	Therefore, \cref{prop:optimalpair:neces:proje} and \cref{prop:optimalpair:neces:o} hold in this case as well.
\end{proof}

\section{Angles between closed  convex cones} \label{sec:anglesclosedconvexcones}

In this section, we characterize the positivity of the minimal angle between
two closed convex cones and study the closedness
of the sum of the two cones.

 \subsection*{Positive angles between two cones}
 \begin{lemma} \label{lemma:cK1K2}
 	Let $K_{1}$ and $K_{2}$ be nonempty closed convex cones in $\mathcal{H}$. Then the following hold:
 	\begin{enumerate}
 		\item \label{lemma:cK1K2:1} If $K_{1} \cap K_{2} \neq \{0\}$, then $c_{0}(K_{1},K_{2}) =1$.
 		\item \label{lemma:cK1K2:0EQ} If $K_{1} \cap K_{2} =\{0\}$,  then $c_{0}(K_{1},K_{2}) =c(K_{1},K_{2}) $.

 		\item \label{lemma:cK1K2:H}  $K_{1} \cap K_{2} =\{0\}$ if and only if $ \overline{K_{1}^{\ominus}+K_{2}^{\ominus}} =\mathcal{H}$.
 	 
 	\end{enumerate}
 \end{lemma}
 
 \begin{proof}
 	\cref{lemma:cK1K2:1}$\&$\cref{lemma:cK1K2:0EQ}: These follow from  \cref{lemma:CD}\cref{lemma:CD:C0=1}  and \cref{fact:coneK1K2}, respectively.

 	\cref{lemma:cK1K2:H}:  By \cref{fact:K1:Km:dualsum}, $K_{1} \cap K_{2} =\{0\} \Leftrightarrow  (K_{1} \cap K_{2})^{\ominus} =\{0\}^{\ominus} =\mathcal{H} \Leftrightarrow
 	\overline{K_{1}^{\ominus}+ K_{2}^{\ominus}} =\mathcal{H}$.
 \end{proof}

\begin{corollary} \label{cor:K1K2c0o}
	Let $K_{1}$ and $K_{2}$ be nonempty closed convex cones in $\mathcal{H}$. Then the following hold:
	\begin{enumerate}	
		\item   \label{cor:K1K2c0o:EQ} Assume that $K_{1} \cap K_{2} =\{0\}$.  Then $c(K_{1},K_{2}) <1 \Leftrightarrow c_{0}(K_{1},K_{2}) <1 \Leftrightarrow $ $K_{1}^{\ominus} +K_{2}^{\ominus}$ is closed.	
			\item  \label{cor:K1K2c0o:closed} Assume that $c_{0}(K_{1},K_{2}) <1$. Then $K_{1}^{\ominus} +K_{2}^{\ominus}$ is closed.
	\end{enumerate}	
\end{corollary}

\begin{proof}
	\cref{cor:K1K2c0o:EQ}: Because $K_{1} \cap K_{2} =\{0\}$ is equivalent to $ (K_{1}\cap K_{2})^{\ominus} =\mathcal{H}$, the desired equivalences follow directly from   \cite[Theorem~2.5]{DH2006II} and \cite[Corollary~4.10]{DH2008III}.
	
	\cref{cor:K1K2c0o:closed}:	 
Note that, by \cref{lemma:cK1K2}\cref{lemma:cK1K2:1}, $c_{0}(K_{1},K_{2}) <1$  	implies $K_{1} \cap K_{2} =\{0\}$. The required result follows from   \cref{cor:K1K2c0o:EQ} above.
\end{proof}

 According  to \cref{example:KMR3} below, we know that under the assumption of \cref{cor:K1K2c0o}, 
the converse statement of \cref{cor:K1K2c0o}\cref{cor:K1K2c0o:closed} doesn't hold even if one of the cones is a linear subspace. 

Moreover, by \cref{example:KMR3}, we see that  the sum of a closed convex cone and  a linear subspace is generally not closed in $\mathbb{R}^{3}$. In addition, the closedness of
$K^{\ominus} + M^{\perp}$  does not imply the closedness of $K + M$, and, by \cref{fact:dualcone}\cref{fact:dualcone:o}$\&$\cref{fact:dualcone:closedconvex},  vice versa. Hence,  \cite[Lemma~2.11]{Deutsch1995} fails when one of the closed linear subspaces is substituted by a closed convex cone.

\begin{example} \label{example:KMR3}
	Suppose $\mathcal{H} =\mathbb{R}^{3}$. Set $K:= \left\{(x_{1}, x_{2},x_{3} ) \in \mathbb{R}^{3} ~:~ \sqrt{x_{1}^{2} + x_{2}^{2}} \leq x_{3} \right\}$ and $M:= \mathbb{R} (1,0,-1)$. Then the following statements hold:
	\begin{enumerate}
		\item \label{example:KMR3:KMo}  $K$ is a closed convex cone and $M$ is a closed linear subspace. Moreover, 
		\begin{align*}
		K^{\ominus}= \left\{ (y_{1}, y_{2},y_{3} ) \in \mathbb{R}^{3} ~:~ y_{3} \leq -\sqrt{y_{1}^{2} + y_{2}^{2}}  \right\} \quad \text{and} \quad  M^{\perp} =  \left\{ (y_{1}, y_{2},y_{3} ) \in \mathbb{R}^{3} ~:~ y_{1} =y_{3}  \right\}.
		\end{align*} 
		\item \label{example:KMR3:notclosed} $K+M$ is not closed.
		\item \label{example:KMR3:c0} $K  \cap M = \mathbb{R}_{+} (-1,0,1) \neq   \{0\}$, $c_{0} (K , M)=1$, and $c  (K , M) =0$.
		
		\item \label{example:KMR3:K+MoH} $K^{\ominus}+M^{\perp} =\{  (x_{1}, x_{2},x_{3} ) \in \mathbb{R}^{3} ~:~ x_{3} \leq x_{1}\}$ is   closed.
		
		\item  \label{example:KMoR3:cap} $K^{\ominus}\cap  M^{\perp} =\{ (y_{1}, y_{2},y_{3} ) \in \mathbb{R}^{3} ~:~ y_{1}=y_{3} \leq 0, y_{2}=0 \} =\mathbb{R}_{+}(-1,0,-1)$ and $ (K^{\ominus} \cap  M^{\perp} )^{\ominus} =\{ (z_{1}, z_{2},z_{3} ) \in \mathbb{R}^{3} ~:~ z_{1} + z_{3} \geq 0 \}$.
		\item  \label{example:KMR3:cKoMo} $c_{0} (K^{\ominus}, M^{\perp})=1$ and $ c (K^{\ominus}, M^{\perp})=0$. 
	\end{enumerate}	
\end{example}

 \begin{proof}
	\cref{example:KMR3:KMo}:  It is clear that   $M$ is a  closed linear subspace, that $K$ is a closed cone and that $M^{\perp} =  \left\{ (y_{1}, y_{2},y_{3} ) \in \mathbb{R}^{3} ~:~ y_{1} =y_{3}  \right\}$.  
	
	Because $K= \left\{(x_{1}, x_{2},x_{3} ) \in \mathbb{R}^{3} ~:~ \sqrt{x_{1}^{2} + x_{2}^{2}} \leq x_{3} \right\} =  \left\{(x_{1}, x_{2},x_{3} ) \in \mathbb{R}^{2} \times \mathbb{R} ~:~ \norm{(x_{1},x_{2})} \leq x_{3} \right\}$,   \cite[Thorem~3.3.6]{BHH1996} implies that  $K$ is convex and  $K^{\ominus} =-K=\left\{ (y_{1}, y_{2},y_{3} ) \in \mathbb{R}^{3} ~:~ y_{3} \leq -\sqrt{y_{1}^{2} + y_{2}^{2}}  \right\}$.

	\cref{example:KMR3:notclosed}: We first show that $(0,1,0) \notin K + M$. Assume to the contrary that  $(0,1,0) \in K + M$. Then there exist $(x_{1}, x_{2},x_{3}) \in K$   and $t \in \mathbb{R}$ such that $(0,1,0)= (x_{1}, x_{2},x_{3})+(t, 0,-t) = (x_{1}+t, x_{2}, x_{3}-t)$. Hence, $t=-x_{1}=x_{3}$ and $x_{2}=1$. Then $x^{2}_{1}+x^{2}_{2} =t^{2}+1>t^{2} = x^{2}_{3}$, which contradicts with the assumption that $(x_{1}, x_{2},x_{3}) \in K $. On the other hand,
	$(0,1,0) =\lim_{t \to \infty} \Big(0, 1+\tfrac{1}{t}, -t+\sqrt{t^{2} +(1+\tfrac{1}{t})^{2}} \Big)  
	= \lim_{t \to \infty} \Big(-t, 1+\tfrac{1}{t}, \sqrt{t^{2} +(1+\tfrac{1}{t})^{2}} \Big) +t(1,0,-1) \in \overline{K +M}$.
		
	Altogether, $K +M$ is not closed.
	
	\cref{example:KMR3:c0}: Because  $(\forall (x_{1}, x_{2},x_{3} ) \in K)$ $x_{3} \geq 0$, it is easy to see that  $K  \cap M = \mathbb{R}_{+} (-1,0,1) \neq   \{0\}$. Then by 	\cref{example:KMR3:KMo} and \cref{lemma:cK1K2}\cref{lemma:cK1K2:1}, $c_{0} (K , M)=1$.	
	In addition,  note that $(K  \cap M )^{\ominus} =( \mathbb{R}_{+} (-1,0,1))^{\ominus} =\{ (y_{1}, y_{2},y_{3})\in \mathbb{R}^{3} ~:~ -y_{1}+y_{3} \leq 0 \}$, and so $K \cap (K  \cap M )^{\ominus} = \mathbb{R}_{+} (1,0,1)  $ and $M \cap (K  \cap M )^{\ominus} = \mathbb{R}_{+} (1,0,-1) $. Hence, $c  (K , M) =c_{0}  (K \cap (K  \cap M )^{\ominus}  , M \cap (K  \cap M )^{\ominus}) =0$.

	\cref{example:KMR3:K+MoH}: Set  $B:= \{  (x_{1}, x_{2},x_{3} ) \in \mathbb{R}^{3} ~:~ x_{3} \leq x_{1}\}$. If $(x_{1}, x_{2},x_{3} ) \in B $, then 
	\begin{align*}
	(x_{1}, x_{2},x_{3} ) = \left(\frac{x_{1}-x_{3}}{2}, 0, \frac{x_{3}-x_{1}}{2}\right) + \left(\frac{x_{1} +x_{3}}{2}, x_{2},\frac{x_{1} +x_{3}}{2}  \right)  \in K^{\ominus}+M^{\perp}, 
	\end{align*}
by \cref{example:KMR3:KMo}. Hence,	  $B \subseteq  K^{\ominus}+M^{\perp}$.
	
	On the other hand,  let $(x_{1}, x_{2},x_{3} ) \in \mathcal{H}\smallsetminus B$, i.e., $x_{3} >x_{1}$. We shall show that $(x_{1}, x_{2},x_{3} ) \notin K^{\ominus}+M^{\perp} $.
	Assume to the contrary that $(x_{1}, x_{2},x_{3} ) \in K^{\ominus}+M^{\perp} $. By  \cref{example:KMR3:KMo}, $ (x_{1}, x_{2},x_{3} ) = (y_{1}, y_{2}, y_{3}   ) + (z_{1}, z_{2},z_{3} )$ where $y_{3} \leq -\sqrt{y_{1}^{2} + y_{2}^{2}}$ and $z_{1}=z_{3}$. Now, $x_{1} = y_{1} +z_{1}$ and $x_{3} =y_{3}+z_{1}$. Combine this with $ x_{3} >x_{1}$ and $y_{3} \leq -\sqrt{y_{1}^{2} + y_{2}^{2}}$ to obtain that $y_{1} \geq -\sqrt{y_{1}^{2} + y_{2}^{2}} \geq 	y_{3} >y_{1}$,
	which is absurd. Altogether, $ K^{\ominus}+M^{\perp}=B$ is closed.

	\cref{example:KMoR3:cap}: Set $C:=\{ (y_{1}, y_{2},y_{3} ) \in \mathbb{R}^{3} ~:~ y_{1}=y_{3} \leq 0, y_{2}=0 \}$. By \cref{example:KMR3:KMo},  $C  \subseteq K^{\ominus} \cap M^{\perp}$. 
	
	Let $(y_{1}, y_{2},y_{3} ) \in K^{\ominus} \cap M^{\perp}$. Then   \cref{example:KMR3:KMo} implies   $(y_{1}, y_{2},y_{3} ) \in K^{\ominus}$ and  $y_{1}=y_{3} \leq 0$. Note that $\abs{y_{3}} \geq \sqrt{y_{1}^{2} + y_{2}^{2}} \geq \abs{y_{1}}$. So $y_{2}=0$ and $K^{\ominus} \cap M^{\perp} =C$. 
	In addition, by the definition of polar cone and $K^{\ominus} \cap M^{\perp} = \{ (y_{1}, y_{2},y_{3} ) \in \mathbb{R}^{3} ~:~ y_{1}=y_{3} \leq 0, y_{2}=0 \}$, we obtain that $ (K^{\ominus} \cap  M^{\perp} )^{\ominus} =\{ (z_{1}, z_{2},z_{3} ) \in \mathbb{R}^{3} ~:~ z_{1} + z_{3} \geq 0 \}$.

	\cref{example:KMR3:cKoMo}: Because, by \cref{example:KMoR3:cap}, $K^{\ominus} \cap M^{\perp} \neq \{0\}   $, by \cref{lemma:cK1K2}\cref{lemma:cK1K2:1},  we have that $c_{0} (K^{\ominus}, M^{\perp})=1$. 	
	Moreover, by \cref{example:KMR3:KMo} and \cref{example:KMoR3:cap}, $K^{\ominus} \cap (K^{\ominus} \cap  M^{\perp} )^{\ominus} =\mathbb{R}_{+} (1,0,-1)  $ and $M^{\perp} \cap (K^{\ominus} \cap  M^{\perp} )^{\ominus} =\mathbb{R}_{+} (1,0,1) $. 
	Hence, by \cref{fact:coneK1K2} and \cref{defn:Angles}, $c (K^{\ominus}, M^{\perp}) = c_{0} \left(  K^{\ominus} \cap (K^{\ominus} \cap  M^{\perp} )^{\ominus} , M^{\perp} \cap (K^{\ominus} \cap  M^{\perp} )^{\ominus}  \right) =0 $.
\end{proof}

The following results imply that in $\mathbb{R}^{n}$, the sufficient conditions in  \cref{lemma:cK1K2}\cref{lemma:cK1K2:1}$\&$\cref{lemma:cK1K2:0EQ} are also necessary conditions.  
 
\begin{proposition} \label{prop:cK1K2}
 Let $\mathcal{K}$ be a finite-dimensional linear subspace of $\mathcal{H}$.  	Let $K_{1}$ and $K_{2}$ be nonempty closed convex cones in $\mathcal{H}$.  Suppose that $K_{1} \subseteq  \mathcal{K}$ or $K_{2} \subseteq  \mathcal{K}$.  Then the following hold:
	\begin{enumerate}
		\item \label{prop:cK1K2:c01} $K_{1} \cap K_{2} \neq \{0\}$ if and only if  $c_{0}(K_{1},K_{2}) =1$.
		\item \label{prop:cK1K2:EQ}   $K_{1} \cap K_{2} =\{0\}$ if and only if  $c_{0}(K_{1},K_{2}) =c(K_{1},K_{2}) $.
		\item \label{prop:cK1K2:c01:EQ} $K_{1} \cap K_{2} = \{0\}$ if and only if  $c_{0}(K_{1},K_{2}) <1$.
	\end{enumerate}
\end{proposition}

\begin{proof}
	\cref{prop:cK1K2:c01}: Assume $c_{0}(K_{1},K_{2}) =1$. Then  $K_{1} \neq \{0\}$ and $K_{2} \neq \{0\}$ and so
 by \Cref{prop:optimalpair:exist,lemma:optimalpair:norm1}, there exist $\bar{x} \in K_{1} \cap \mathbf{S}_{\mathcal{H}}  $ and $\bar{y} \in K_{2}  \cap \mathbf{S}_{\mathcal{H}}   $ such that $\innp{\bar{x} ,\bar{y}  } = c_{0}(K_{1},K_{2}) =1$.
Hence, $\bar{x}\neq 0$, $\bar{y}\neq 0$ and 
	\begin{align*}
	\norm{\bar{x}- \bar{y} }^{2} =\norm{\bar{x}}^{2} -2 	\innp{\bar{x}, \bar{y}}  + \norm{\bar{y}}^{2}  \leq 1-2+1=0,
	\end{align*}
which implies that $	\bar{x}= \bar{y} \in K_{1}\cap K_{2} \smallsetminus\{0\}$.
Therefore,  $K_{1} \cap K_{2} \neq \{0\}$.
	
	Moreover, the reverse direction holds by \cref{lemma:cK1K2}\cref{lemma:cK1K2:1}. Hence,  \cref{prop:cK1K2:c01} holds.
	
	\cref{prop:cK1K2:EQ}: Suppose that $c_{0}(K_{1},K_{2}) =c(K_{1},K_{2}) $. Assume to the contrary that $K_{1} \cap K_{2} \neq \{0\}$. Then by \cref{prop:cK1K2:c01} above, $c(K_{1},K_{2})=c_{0}(K_{1},K_{2}) =1 $, which, by  \cref{fact:coneK1K2},  implies that 
	\begin{align*}
	1 =c(K_{1},K_{2}) = c_{0} \left(  K_{1}  \cap (K_{1}\cap K_{2})^{\ominus}, K_{2} \cap (K_{1}\cap K_{2})^{\ominus} \right). 
	\end{align*}
	Note that by \cref{fact:dualcone}\cref{fact:dualcone:o}, $K_{1}  \cap (K_{1}\cap K_{2})^{\ominus}$ and $K_{2}  \cap (K_{1}\cap K_{2})^{\ominus}$ are nonempty closed convex cones. Apply \cref{prop:cK1K2:c01} above with $K_{1}=K_{1}  \cap (K_{1}\cap K_{2})^{\ominus}$ and $K_{2} = K_{2} \cap (K_{1}\cap K_{2})^{\ominus}$ to obtain that 
	\begin{align*}
	\{0\} = (K_{1}\cap K_{2})  \cap (K_{1}\cap K_{2})^{\ominus}=\left( K_{1}  \cap (K_{1}\cap K_{2})^{\ominus} \right)\cap  \left( K_{2}  \cap (K_{1}\cap K_{2})^{\ominus} \right) \neq \{0\},
	\end{align*}
	which is absurd. 
	Hence, $K_{1} \cap K_{2} = \{0\}$. In addition, the converse direction holds by  \cref{lemma:cK1K2}\cref{lemma:cK1K2:0EQ}.  
	
	\cref{prop:cK1K2:c01:EQ}:  Note that by \cref{fact:AnglesProperties}\cref{fact:AnglesProperties:between0and1}, $c_{0}(K_{1},K_{2}) \neq 1 \Leftrightarrow c_{0}(K_{1},K_{2}) <1$. So, the desired result is from \cref{prop:cK1K2:c01}.
\end{proof}

The following result yields that the cosine of the angle defined in \cref{defn:Angles} between two closed convex cones is always strictly less than $1$ in a finite-dimensional space. 
\begin{proposition} \label{prop:c<1}
	Let $\mathcal{K}$ be a finite-dimensional linear subspace of $\mathcal{H}$.  	Let $K_{1}$ and $K_{2}$ be nonempty closed convex cones of $\mathcal{H}$.  Suppose that $K_{1} \subseteq  \mathcal{K}$ or $K_{2} \subseteq  \mathcal{K}$.  Then $	c(K_{1},K_{2}) <1$.
\end{proposition}

\begin{proof}
	If $K_{1} \cap K_{2} \neq \{0\}$, then by \cref{prop:cK1K2}\cref{prop:cK1K2:EQ} and \cref{lemma:CD}\cref{lemma:CD:c0:c}, we have that $c(K_{1},K_{2}) <c_{0}(K_{1},K_{2})  \leq 1$. 	
Assume that $K_{1} \cap K_{2} = \{0\}$.	Then combine \cref{lemma:CD}\cref{lemma:CD:c0:c} with \cref{prop:cK1K2}\cref{prop:cK1K2:c01:EQ} to see that $c(K_{1},K_{2}) \leq c_{0}(K_{1},K_{2}) <1$. 
\end{proof}

In the following result, we characterize when the cosine of the minimal angle between two nonempty closed convex cones is strictly less than $1$.

\begin{theorem} \label{theorem:c0K1K2EQ}
 	Let $K_{1}$ and $K_{2}$ be nonempty closed convex cones in $\mathcal{H}$ such that $K_{1} \neq \{0\}$ and $K_{2} \neq \{0\}$. Then the following statements are equivalent:
	\begin{enumerate}
		\item \label{theorem:c0K1K2EQ:c0} $c_{0}(K_{1},K_{2}) <1$.
		\item \label{theorem:c0K1K2EQ:norm} $\sup \left\{    \innp{\Pro_{K_{1}}x,\Pro_{K_{2}}\Pro_{K_{1}}x} ~:~  x \in  \mathbf{S}_{\mathcal{H}} \right\} <1$.
		\item \label{theorem:c0K1K2EQ:d} $\dist \left(  K_{1} \cap \mathbf{S}_{\mathcal{H}}, K_{2} \cap \mathbf{S}_{\mathcal{H}} \right) >0$.
		\item \label{theorem:c0K1K2EQ:gamma} $\inf \left\{   \frac{ \norm{x-y}} {\norm{x} +\norm{y}}~:~ x \in K_{1}, ~ y \in K_{2}, ~(x,y) \neq (0,0) \right\} >0$.
	\end{enumerate}
\end{theorem}

\begin{proof}
According to	\cref{theorem:c0}\cref{theorem:c0:P}, we have 	\cref{theorem:c0K1K2EQ:c0} $\Leftrightarrow$ \cref{theorem:c0K1K2EQ:norm}. Moreover, recalling \cref{fact:coBS} and applying \cref{lemma:ineq} with $C$ and $D$ replaced by $K_{1}$ and $K_{2}$ respectively, we obtain that \cref{theorem:c0K1K2EQ:c0} $\Leftrightarrow$ \cref{theorem:c0K1K2EQ:d}  $\Leftrightarrow$ \cref{theorem:c0K1K2EQ:gamma}.	
\end{proof}

\subsection*{Closedness of sums of cones}

Consider the closed convex cone $K:=\left\{(x_{1}, x_{2},x_{3} ) \in \mathbb{R}^{3} ~:~ \sqrt{x_{1}^{2} + x_{2}^{2}} \leq x_{3} \right\}  $ and the closed linear subspace $M:= \mathbb{R} (1,0,-1)$.  By \cref{example:KMR3}\cref{example:KMR3:notclosed}, we know that $K+M$ is not closed. We provide sufficient conditions for the closedness of the sum of two closed convex cones below. In view of \cref{example:KMR3}, $K^{\ominus}- M^{\perp}$ is closed but $c_{0}(K^{\ominus}, M^{\perp}) =1$, which implies that the converse of  \cref{theorem:K1K2closed} below	 is generally not true.

The   idea of the proof of \cref{theorem:K1K2closed}   is from that of  \cite[Theorem~2.12(1)$\Rightarrow$(2)]{Deutsch1995} which is on linear subspaces.

\begin{theorem} \label{theorem:K1K2closed}	
		Let $K_{1}$ and $K_{2}$ be  nonempty closed convex cones in $\mathcal{H}$. Assume that $c_{0}(K_{1}, K_{2}) <1$. Then $K_{1} - K_{2} $ is closed. 
\end{theorem}

\begin{proof}
	Take $\bar{z} \in \overline{K_{1} -K_{2}} $. Then  there exist  sequences   $(x_{k})_{k \in \mathbb{N}}$ in $K_{1}  $ and $ (y_{k})_{k \in \mathbb{N}} $ in $K_{2} $ such that $  x_{k}-y_{k} \to \bar{z}$.
By \cref{fact:AnglesProperties}\cref{fact:AnglesProperties:ineq}, we know that $(\forall x \in K_{1} ) (\forall y \in K_{2} ) \quad \innp{x,y} \leq c_{0}(K_{1}, K_{2})  \norm{x} \norm{y}$. So,
	$(\forall k \in \mathbb{N})$  $ \norm{x_{k}-y_{k}}^{2} = \norm{x_{k}}^{2} +\norm{y_{k}}^{2} - 2 \innp{x_{k},y_{k}}  \geq (\norm{x_{k}} -\norm{y_{k}} )^{2}+ 2(1- c_{0}(K_{1},K_{2})) \norm{x_{k}} \norm{y_{k}}$,
	which, combining with  the boundedness of the convergent sequence $(x_{k}-y_{k})_{k \in \mathbb{N}}$ and the assumption $c_{0}(K_{1}, K_{2}) <1$, yields that the sequences $  (\norm{x_{k}} -\norm{y_{k}} )_{k \in \mathbb{N}}$ and $( \norm{x_{k}} \norm{y_{k}})_{k \in \mathbb{N}}$ are bounded. 
	Note that $(\forall k \in \mathbb{N})$ $(\norm{x_{k}} -\norm{y_{k}} )^{2}+ 2\norm{x_{k}} \norm{y_{k}} =\norm{x_{k}}^{2} +\norm{y_{k}}^{2}  $.  Hence,  $(x_{k})_{k \in \mathbb{N}}$ and  $(y_{k})_{k \in \mathbb{N}}$ are bounded. 
	
	Because $K_{1}$ and $K_{2}$ are nonempty closed and convex sets, by \cite[Lemma~2.45 and Corollary~3.35]{BC2017} and by passing to a subsequence if necessary, there exist $\bar{x} \in K_{1}$ and $\bar{y} \in K_{2}$ such that $x_{k} \weakly \bar{x}$ and $y_{k} \weakly \bar{y}$. Combine this with the result that $ x_{k}-y_{k} \to \bar{z}$ to obtain that $\bar{z} = \bar{x}-\bar{y} \in K_{1}-K_{2}$. Hence, $ \overline{K_{1} -K_{2}}  \subseteq K_{1} -K_{2}$. Therefore,    $K_{1}-K_{2}$ is closed.  
\end{proof}
\begin{corollary}\label{cor:K1K2closed}	
 	Let $K_{1}$ and $K_{2}$ be  nonempty closed convex cones in $\mathcal{H}$.  Assume that one of the following items holds:
\begin{enumerate}
		\item  \label{theorem:K1K2closed:c0}  $c_{0}(K_{1}, -K_{2}) <1$.  
		\item \label{theorem:K1K2closed:sqrt}  $\sup \left\{   \innp{\Pro_{K_{1}}x,\Pro_{(-K_{2})}\Pro_{K_{1}}x} ~:~  x \in  \mathbf{S}_{\mathcal{H}} \right\} <1$.
	\item  \label{theorem:K1K2closed:d}  $\dist \left(  K_{1} \cap \mathbf{S}_{\mathcal{H}}, -K_{2} \cap \mathbf{S}_{\mathcal{H}} \right) >0  $.
	\item \label{theorem:K1K2closed:gamma} $\inf \left\{   \frac{\norm{x+y}}{\norm{x} +\norm{y} } ~:~ x \in K_{1}, y \in K_{2}, (x,y) \neq (0,0)    \right\} >0$.
	\item \label{theorem:K1K2closed:cap0}	
Let $\mathcal{K}$ be a finite-dimensional linear subspace of $\mathcal{H}$, and suppose that $K_{1} \subseteq  \mathcal{K}$ or $K_{2} \subseteq  \mathcal{K}$ and that $K_{1} \cap (-K_{2}) =\{0\}$.
	\end{enumerate}
Then $K_{1} + K_{2} $ is closed. 
\end{corollary} 

\begin{proof}
	If $K_{1}=\{0\}$ or $K_{2} =\{0\}$, then clearly, $K_{1} + K_{2} $ is closed. Assume that $K_{1} \neq \{0\}$ and $K_{2} \neq \{0\}$.
	Note that $K_{2}$ is a closed convex cone if and only if  $-K_{2}$ is a closed convex cone. 
By \cref{theorem:c0K1K2EQ}, we know that the conditions 
	\cref{theorem:K1K2closed:c0},  \cref{theorem:K1K2closed:sqrt},  \cref{theorem:K1K2closed:d} and  \cref{theorem:K1K2closed:gamma} are equivalent. Moreover, by \cref{prop:cK1K2}\cref{prop:cK1K2:c01:EQ},
	the condition  \cref{theorem:K1K2closed:cap0} implies   \cref{theorem:K1K2closed:c0}. Combine these results with \cref{theorem:K1K2closed}	 to complete the proof. 
\end{proof}

\begin{remark}
	Consider \cref{cor:K1K2closed}. 
	\begin{enumerate}
		\item The fact that \cref{theorem:K1K2closed:d}  implies the closedness of $K_{1} + K_{2} $ follows also from  \cite[Theorem~3.2]{Beutner2007}.
		\item The result that the conditions \cref{theorem:K1K2closed:gamma}  and \cref{theorem:K1K2closed:cap0}  imply the closedness of $K_{1} + K_{2} $ is the second part of \cite[Proposition~4.1]{SS2016I} which is in a  Euclidean space.
	\end{enumerate}
In conclusion, by using the tool of minimal angle, we deduce the sufficient conditions for the closedness of the sum of two closed convex cones in \cite{Beutner2007} and \cite{SS2016I}. 
\end{remark}

\subsection*{Intersections of cones}
	To prove the main result \cref{theorem:KominusMperpNeq0} in this subsection, we need the following easy result.
	\begin{lemma} \label{lemma:Ku}
		Let $K$ be a  convex subset of $\mathcal{H}$ and let $u \in \mathcal{H}$. Assume that there exist $x \in K$ and $y \in K$ such that $\innp{x,u}>0$ and $\innp{y,u} <0$.  Then there exists $t \in \left]0,1\right[$ such that  $tx+(1-t)y \in K \cap \{u  \}^{\perp}$. 
	\end{lemma}
\begin{proof}
 Set  $f: \left[ 0,1\right] \to \mathbb{R}: t \mapsto \innp{y +t(x-y),u}$. Then the required result follows easily from the intermediate value theorem applied to $f$.
\end{proof}

\begin{theorem} \label{theorem:KominusMperpNeq0}
	Let $K_{1} $ and $K_{2}$  be  nonempty closed convex cones in $\mathcal{H}$.  Assume that $K_{1}$ is not a linear subspace and that $K_{1} \cap K_{2}=\{0\}$. Assume that one of the following items hold:
	\begin{enumerate}
		\item \label{theorem:KominusMperpNeq0:M} There exists $u \in \mathcal{H}$ such that $K_{2} =\{u  \}^{\perp}$.
		\item \label{theorem:KominusMperpNeq0:H} There exists $u \in \mathcal{H}$ such that $ H:=\{u  \}^{\perp}$  satisfies $K_{2} \subseteq H$ and $H\cap K_{1} =\{0\}$.
		\item \label{theorem:KominusMperpNeq0:K} Let $\mathcal{K}$ be a finite-dimensional linear subspace of $\mathcal{H}$.  Suppose that $K_{1} \subseteq  \mathcal{K}$ or $K_{2} \subseteq  \mathcal{K}$.
	\end{enumerate}
	Then 
	\begin{align*}
	K^{\ominus}_{1} \cap K^{\oplus}_{2} \neq \{0\} \quad \text{and} \quad K^{\oplus}_{1} \cap K^{\ominus}_{2} \neq \{0\}.
	\end{align*}
\end{theorem}

\begin{proof}
	Because, by \cref{lemma:polardual}\cref{lemma:polardual:-C},  $K^{\oplus}_{1} \cap K^{\ominus}_{2} =- (K^{\ominus}_{1} \cap K^{\oplus}_{2} ) $, we only need to prove that  $K^{\ominus}_{1} \cap K^{\oplus}_{2} \neq \{0\}$. 
	
	If $K_{2} =\{0\}$, then $K_{2}^{\oplus} =\mathcal{H}$. Because $K_{1}$ is not a linear subspace,  we know that  $K_{1} \neq \mathcal{H}$ and   that $K_{1}^{\ominus} \cap K_{2}^{\oplus} = K_{1}^{\ominus} \neq \{0\} $. Hence, in the rest of the proof, we assume that $K_{2} \neq \{0\}$.

	Assume	\cref{theorem:KominusMperpNeq0:M} holds: We separate the proof into the following two steps:

	\emph{Step~1:} In this part, we show that 
	\begin{align} \label{theorem:KominusMperpNeq0:M:eitheror}
 (\forall x \in K_{1}) ~ \innp{u,x} \leq 0 \quad \text{or} \quad (\forall x \in K_{1})  ~ \innp{u,x} \geq 0.
	\end{align}
	Assume to the contrary that there exists $x_{1} \in K_{1}$ and $x_{2} \in K_{1}$ such that $\innp{x_{1},u}>0$ and $\innp{x_{2},u} <0$. Then apply \cref{lemma:Ku} with $K=K_{1}$ to see that there exists $\bar{t} \in \left]0,1\right[$ such that $\bar{t} x_{1}+(1-\bar{t}) x_{2} \in K_{1} \cap K_{2}$.
	If $\bar{t} x_{1}+(1-\bar{t}) x_{2} \neq 0$, then $\bar{t} x_{1}+(1-\bar{t}) x_{1} \in (K_{1} \cap K_{2}) \smallsetminus \{0\}$, which contradicts   the assumption that $K_{1} \cap K_{2} =\{0\}$. So, we know that \cref{theorem:KominusMperpNeq0:M:eitheror} is true.
	
	Now assume $ \bar{t} x_{1}+(1-\bar{t}) x_{2}  =0$. Then  $x_{1} = - \frac{1-\bar{t} }{\bar{t} }  x_{2} \in -K_{1}$,
	which implies that $x_{1}  \in K_{1} \cap ( -K_{1})$. Combine this with the assumption that  $K_{1} $ is a closed convex cone to see that $\spn \{x_{1}\}=\mathbb{R} \cdot x_{1}   \subseteq K_{1}$. 
	
	Because $K_{1}$ is not a linear subspace,  we see that $ K_{1} \smallsetminus\spn \{x_{1}\}  \neq \varnothing$.
	Take $z \in K_{1} \smallsetminus \spn \{x_{1}\}$.  Because $K_{1} \cap K_{2} =\{0\}$ and $z \neq 0$, we have either $\innp{z,u} >0$ or $\innp{z,u} <0$.
	
	If $\innp{z,u} >0$, then applying \cref{lemma:Ku} with $K=K_{1}$, $x =z$ and $y=-x_{1}$, we get that there exists $\tilde{t}   \in \left]0,1\right[$ such that $\tilde{t} z -(1-\tilde{t}) x_{1} \in K_{1} \cap K_{2}$.  If $\tilde{t} z -(1-\tilde{t}) x_{1}=0$, then $z =\frac{1-\tilde{t} }{\tilde{t} } x_{1} \in  \spn \{x_{1}\}$, which contradicts    that $z \in K_{1} \smallsetminus \spn \{x_{1}\}$. Hence, $\tilde{t} z -(1-\tilde{t}) x_{1} \in (K_{1} \cap K_{2}) \smallsetminus \{0\}$, which contradicts   the assumption that $K_{1} \cap K_{2} =\{0\}$. Hence,  \cref{theorem:KominusMperpNeq0:M:eitheror} is true in this case.

If $\innp{z,u} <0$, then an analogous argument yields a contradiction.
Hence,  in this case   \cref{theorem:KominusMperpNeq0:M:eitheror} holds  as well.

	Altogether, in all cases,  \cref{theorem:KominusMperpNeq0:M:eitheror} holds.
	
	\emph{Step~2:} 
	Note that  if $u=0$, then using  the assumption that $K_{1}$ is not a linear subspace, we have that  $K_{1} \cap K_{2} =K_{1} \neq \{0\}$, which  contradicts   the assumption $K_{1} \cap K_{2} =\{0\}$. Hence, $u \neq 0$.  Moreover, because $K_{2} =\{u  \}^{\perp}$, by \cref{lemma:polardual}\cref{lemma:polardual:llinear}, $K^{\oplus}_{2}=K^{\perp}_{2}=\spn \{u\}$.
Therefore, by \cref{theorem:KominusMperpNeq0:M:eitheror}, we have exactly the following two cases:
	
	\emph{Case~1}: $(\forall x \in K_{1}) $  $\innp{u,x} \leq 0$. Then $u \in (K_{1}^{\ominus} \cap K_{2}^{\perp}) \smallsetminus \{0\}$.
	
	\emph{Case~2}:  $(\forall x \in K_{1}) $  $\innp{u,x} \geq 0$. Then $-u \in (K_{1}^{\ominus} \cap K_{2}^{\perp}) \smallsetminus \{0\}$. 
	
	Altogether, $K_{1}^{\ominus} \cap K_{2}^{\perp} \neq \{0\} $.
	
	Assume  \cref{theorem:KominusMperpNeq0:H} holds:  Because $K_{2} \subseteq H$ implies that $H^{\perp} =H^{\oplus} \subseteq K_{2}^{\oplus} $, apply \cref{theorem:KominusMperpNeq0:M} with $K_{2} =H$ to obtain that $	\{0\} \neq K_{1}^{\ominus} \cap H^{\perp} \subseteq K_{1}^{\ominus} \cap K_{2}^{\oplus}$,
	which implies that $K_{1}^{\ominus} \cap K_{2}^{\oplus} \neq \{0\}$.
	
	Assume   \cref{theorem:KominusMperpNeq0:K} holds: Assume to the contrary that $K_{1}^{\ominus} \cap K_{2}^{\oplus} = \{0\}$. Then by \cref{lemma:cK1K2}\cref{lemma:cK1K2:H}, \cref{fact:dualcone}\cref{fact:dualcone:closedconvex} and \cref{lemma:polardual}\cref{lemma:polardual:-C}, we know that 
	\begin{align}\label{eq:theorem:KominusMperpNeq0:K:H}
	\mathcal{H} = \overline{K_{1}^{\ominus \ominus} + K_{2}^{\oplus\ominus}} =\overline{K_{1}-K_{2}}.
	\end{align}
	Because $K_{1}\cap K_{2} =\{0\}$, by \cref{prop:cK1K2}\cref{prop:cK1K2:c01:EQ} and \cref{theorem:K1K2closed}, we have that 	
	\begin{align*}
	K_{1}\cap K_{2} =\{0\} \Leftrightarrow c_{0}(K_{1},K_{2})<1 \Rightarrow K_{1}-K_{2} = \overline{K_{1}-K_{2} },
	\end{align*}
	which, combining with \cref{eq:theorem:KominusMperpNeq0:K:H}, implies that 
	\begin{align}\label{eq:theorem:KominusMperpNeq0:K:K+M}
	\mathcal{H} =K_{1} -K_{2}.
	\end{align}
	Because $K_{1}$ is a nonempty convex cone but not a linear subspace, by \cref{lemma:K}, $-K_{1} \not \subseteq K_{1}$. Hence,  there exists $\bar{y} \in K_{1}$ such that $- \bar{y} \notin K_{1}$. Take $\bar{x} \in K_{2} \smallsetminus \{0\}$. Then  $	\bar{x} - \bar{y} \in \mathcal{H} \stackrel{\cref{eq:theorem:KominusMperpNeq0:K:K+M}}{=} K_{1}-K_{2}, $
	which implies that there exist $\tilde{y} \in K_{1}$ and  $\tilde{x} \in K_{2}$  such that 
	\begin{align}\label{eq:theorem:KominusMperpNeq0:eq}
	\bar{x} - \bar{y} =   \tilde{y} -\tilde{x}.
	\end{align}
	Because $K_{1}$ and $K_{2}$ are nonempty closed convex cones, by \cref{lemma:K}, $K_{1}+K_{1}=K_{1}$ and $K_{2}+K_{2}=K_{2}$. Hence, $\bar{x} + \tilde{x} \stackrel{\cref{eq:theorem:KominusMperpNeq0:eq}}{=} \tilde{y} +\bar{y} \in K_{2} \cap K_{1}$.
	If $\bar{x} + \tilde{x} =0$, then by \cref{eq:theorem:KominusMperpNeq0:eq}, $- \bar{y} =   \tilde{y}-(\bar{x} + \tilde{x}) =\tilde{y} \in K_{1}$, which contradicts   $ - \bar{y} \notin K_{1}$. Hence, $\bar{x} + \tilde{x}  \in (K_{1} \cap K_{2}) \smallsetminus \{0\}$, which contradicts   the assumption that $K_{1} \cap K_{2}=\{0\}$. Therefore, under the condition \cref{theorem:KominusMperpNeq0:K}, we have also that  $K_{1}^{\ominus} \cap K_{2}^{\oplus} \neq \{0\}$.
\end{proof}

As an application of \cref{theorem:KominusMperpNeq0}, in the following \cref{coro:KominusMperpNeq0}, we show that under the assumptions of \cref{theorem:KominusMperpNeq0}, $c_{0}(K_{1},K_{2})<1$ and $c_{0}(K_{1}^{\ominus}, K_{2}^{\oplus}) <1$ cannot occur together. 
Hence, under the assumptions of \cref{coro:KominusMperpNeq0}, $c_{0} (K_{1}, K_{2}) <1$ implies that $c_{0}(K_{1}^{\ominus}, K_{2}^{\oplus}) =1 > c_{0} (K_{1}, K_{2})$, which reduces to \cite[Lemma~2.14]{Deutsch1995} when $M$ and $N$ are cones and $X=\mathcal{H}$.
\begin{corollary}  \label{coro:KominusMperpNeq0}
	Let $K_{1}$ and $K_{2}$ be  nonempty closed convex cones in $\mathcal{H}$.  Suppose  that $K_{1}$ is not a linear subspace. Assume that one of the following items holds:
	\begin{enumerate}[label=(\alph*)]
		\item \label{coro:KominusMperpNeq0:M} There exists $u \in \mathcal{H}$ such that $K_{2}=\{u  \}^{\perp}$.
		\item \label{coro:KominusMperpNeq0:H} There exists $u \in \mathcal{H}$ such that $K_{2} \subseteq H:=\{u  \}^{\perp}$ and that $H\cap K_{1} =\{0\}$.
		\item \label{coro:KominusMperpNeq0:K} Let $\mathcal{K}$ be a finite-dimensional linear subspace of $\mathcal{H}$.  Suppose that $K_{1} \subseteq  \mathcal{K}$ or $K_{2} \subseteq  \mathcal{K}$.
	\end{enumerate}
	Then the following hold:
	\begin{enumerate}
		\item \label{coro:KominusMperpNeq0:dichotomy} If $K_{1} \cap K_{2} \neq \{0\}$, then $c_{0}(K_{1},K_{2})=1$; if $K_{1} \cap K_{2}=\{0\}$, then $c_{0}(K_{1}^{\ominus}, K_{2}^{\oplus}) =1$ and $c_{0}(K_{1}^{\oplus}, K_{2}^{\ominus}) =1$. 
		\item \label{coro:KominusMperpNeq0:c0} If $c_{0}(K_{1}, K_{2}) <1$, then $c_{0}(K_{1}^{\ominus}, K_{2}^{\oplus}) =1$ and $c_{0}(K_{1}^{\oplus}, K_{2}^{\ominus}) =1$.
	\end{enumerate} 
\end{corollary}

\begin{proof}
	\cref{coro:KominusMperpNeq0:dichotomy}:	If $K_{1} \cap K_{2} \neq \{0\}$, then by \cref{lemma:cK1K2}\cref{lemma:cK1K2:1}, $c_{0}(K_{1}, K_{2})=1$.
	
	Assume $K_{1} \cap K_{2}=\{0\}$. 
	By the assumptions and  \cref{theorem:KominusMperpNeq0}, we know that $K_{1}^{\ominus} \cap K_{2}^{\oplus} \neq \{0\}$ and $K^{\oplus}_{1} \cap K^{\ominus}_{2} \neq \{0\}$.
	Combine this with  \cref{lemma:cK1K2}\cref{lemma:cK1K2:1} to obtain $c_{0}(K_{1}^{\ominus}, K_{2}^{\oplus}) =1 $ and $c_{0}(K_{1}^{\oplus}, K_{2}^{\ominus}) =1$. 
	
	\cref{coro:KominusMperpNeq0:c0}:  This is immediate from    \cref{coro:KominusMperpNeq0:dichotomy}.
\end{proof}

According to related definitions and \cref{lemma:CD}\cref{lemma:CD:--}, it is easy to deduce   two examples below. 
\begin{example}\label{exam:K1K2}
	Suppose that $\mathcal{H} =\mathbb{R}^{2}$. Set $K_{1}:=\mathbb{R}^{2}_{+}$ and $K_{2}:= \{ (x_{1}, x_{2}) \in \mathbb{R}^{2} ~:~ -x_{1} \geq x_{2}   \}$. Then the following statements hold:
	\begin{enumerate}
		\item \label{exam:K1K2:ominus} $K_{1}^{\ominus}= \mathbb{R}^{2}_{-}$, $K_{2}^{\ominus}=\mathbb{R}_{+} (1,1)$, and $K_{2}^{\oplus} =\mathbb{R}_{+}(-1,-1)$.
		\item  \label{exam:K1K2:cap} $K_{1} \cap K_{2} =\{0\}$,  $K_{1}^{\ominus} \cap K_{2}^{\ominus} =\{0\}$, $K_{1}^{\ominus} \cap K_{2}^{\oplus} =K_{2}^{\oplus}$, $K_{1} +K_{2} =\mathbb{R}^{2}$, and $K_{1} - K_{2} =- K_{2}$.
		\item \label{exam:K1K2:ineq} $c  (K_{1}, K_{2}) =c_{0} (K_{1}, K_{2}) =\frac{\sqrt{2}}{2} > 0=c_{0}  (K_{1}^{\ominus}, K_{2}^{\ominus}) =c_{0}  (K_{1}^{\oplus}, K_{2}^{\oplus}) =c (K_{1}^{\ominus}, K_{2}^{\ominus}) =c  (K_{1}^{\oplus}, K_{2}^{\oplus})$.
		
		\item \label{exam:K1K2:ineqoplus} $c_{0} (K_{1}^{\ominus},K_{2}^{\oplus}) =1$, $c  (K_{1}^{\ominus},K_{2}^{\oplus}) =0$, $c_{0} (K_{1},K_{2}) <  c_{0} (K_{1}^{\ominus},K_{2}^{\oplus}) $, and  $c  (K_{1},K_{2}) > c  (K_{1}^{\ominus},K_{2}^{\oplus}) $.
	\end{enumerate}
\end{example}

\begin{example} \label{exam:11NEQ}
	Suppose $\mathcal{H} =\mathbb{R}^{2}$. Set $K:= \left\{(x_{1}, x_{2} ) \in \mathbb{R}^{2} ~:~  x_{1} \geq  x_{2}  \geq 0 \right\}$ and $M:= \mathbb{R} (1,0)$. Then the following statements hold:
	\begin{enumerate}
		\item $K^{\ominus} = \{ (x_{1},x_{2}) ~:~ -x_{1} \geq x_{2} \text{ and } x_{1} \leq 0 \} $, and $  M^{\perp}= \mathbb{R}(0,1) $.
		\item  $K \cap M=\mathbb{R}_{+}(1,0) \neq \{0\}$,  $K+M =\{ (x_{1},x_{2}) ~:~ x_{2} \geq 0 \} $, and $K^{\ominus} \cap M^{\perp}  =\mathbb{R}_{+}(0,-1) \neq \{0\}$.
		\item $c_{0} (K,M) =1$,  $c  (K,M) =0$, $ c_{0} (K^{\ominus} , M^{\perp} ) =1$, and  $c  (K^{\ominus} , M^{\perp} ) =\frac{\sqrt{2}}{2}$.
	\end{enumerate}
\end{example}

\begin{remark}
	\begin{enumerate}	
		\item Let $K_{1} $ and $K_{2}$  be  nonempty closed convex cones in $\mathbb{R}^{2}$. 	According to \cref{exam:K1K2}\cref{exam:K1K2:cap},  even if both  $K_{1} $ and $K_{2}$  are not linear subspaces, $K_{1} \cap  K_{2} =\{0\}$ and  $ K_{1}^{\ominus} \cap  K_{2}^{\ominus} =\{0\}$  appear simultaneously. 
		Because  $K_{1}^{\oplus} \cap  K_{2}^{\oplus} = -(K_{1}^{\ominus} \cap  K_{2}^{\ominus})$,  we know that $K_{1} \cap  K_{2} =\{0\}$ and $K_{1}^{\oplus} \cap  K_{2}^{\oplus} =\{0\}$ occur together in \cref{exam:K1K2} as well.  Therefore, we conclude that in \cref{theorem:KominusMperpNeq0}, $K^{\ominus}_{1} \cap K^{\oplus}_{2} \neq \{0\}$ can not be replaced by $ K_{1}^{\ominus} \cap  K_{2}^{\ominus} \neq \{0\}$ or $K_{1}^{\oplus} \cap  K_{2}^{\oplus} \neq \{0\}$.
		\item By  \cref{prop:cK1K2}\cref{prop:cK1K2:c01:EQ} and \cref{theorem:c0K1K2EQ}, the conditions 
		\cref{theorem:K1K2closed:c0}, \cref{theorem:K1K2closed:sqrt}, \cref{theorem:K1K2closed:d}, \cref{theorem:K1K2closed:gamma} and \cref{theorem:K1K2closed:cap0} in \cref{cor:K1K2closed}	 are equivalent in finite-dimensional spaces. Note that in \cref{exam:11NEQ}, $K+M$ is closed, but $c_{0}(K,M) =1$ and $K \cap M \neq \{0\}$. Hence, 
		$K+M$ is closed is generally not a sufficient condition for any of the conditions \cref{theorem:K1K2closed:c0}, \cref{theorem:K1K2closed:sqrt}, \cref{theorem:K1K2closed:d}, \cref{theorem:K1K2closed:gamma} and \cref{theorem:K1K2closed:cap0} in \cref{cor:K1K2closed}, even if one of the underlying two cones is a  linear subspace.
		
		\item  Suppose that the $K_{1}$ and $K_{2}$ in \cref{theorem:KominusMperpNeq0} are respectively $K$ and $M$ in \cref{exam:11NEQ}. Now, we have that \cref{theorem:KominusMperpNeq0:M}, \cref{theorem:KominusMperpNeq0:H} and \cref{theorem:KominusMperpNeq0:K} in \cref{theorem:KominusMperpNeq0} hold,  and that $K^{\ominus} \cap M^{\perp}    \neq \{0\}$, but  $K \cap M  \neq \{0\}$. Hence, we conclude that under conditions \cref{theorem:KominusMperpNeq0:M}, \cref{theorem:KominusMperpNeq0:H} and \cref{theorem:KominusMperpNeq0:K} in \cref{theorem:KominusMperpNeq0},  $K^{\ominus}_{1} \cap K^{\oplus}_{2} \neq \{0\}$ generally does not imply $K_{1} \cap K_{2}=\{0\}$, although the converse statement was shown in \cref{theorem:KominusMperpNeq0}. 
	\end{enumerate}
\end{remark}

\section{Applications and conclusion}\label{sec:applications}
One of the most important applications of minimal angles is in the study of  convergence rates of cyclic projections algorithms.
Deutsch and Hundal  showed that the rate of convergence for the cyclic projections algorithm onto an intersection of finitely many closed convex sets can be
described by the  \enquote{norm} of the composition of projectors onto certain
sets constructed by these convex sets, and that under some conditions,  the rate of convergence is the square of  the cosine of angle   between related closed convex cones. To precisely state their result, we require two definitions: Let $F:\mathcal{H} \to \mathcal{H}$, let $A$ be a nonempty closed convex subset of $\mathcal{H}$ and let $\epsilon \in \mathbb{R}_{+}$. Then $\norm{F}:=\sup \{ \frac{\norm{F(x)}}{\norm{x}}  ~:~ x\in \mathcal{H} \smallsetminus \{0  \} \}$, and the  \emph{$\epsilon$-dual cone of $A$} is the set $A^{\ominus,\epsilon} := \cone \{ x-\Pro_{A}x ~:~ x \in \mathbf{B}[0 ; \epsilon]  \}$. We are now ready to state Deutsch and Hundal's result.
\begin{fact} {\rm \cite[Theorem~5.21 and Corollary~5.22]{DH2006II}}
Let  $C$ and $D$ be closed convex sets with $C \cap D \neq \varnothing$ and let $ x_{0} \in \mathcal{H}$. Set	$(\forall k \in \mathbb{N})$ $x_{k+1} := \Pro_{D}\Pro_{C} x_{k}$. Then there exists $\bar{x} \in C \cap D$ such that $x_{k} \weakly \bar{x}$,  and for any $\epsilon \in \mathbb{R}_{++}$ with $\epsilon \geq \norm{x_{0} -\bar{x}}$, 
	\begin{align*}
	(\forall k \in \mathbb{N})	\quad	\norm{x_{k+1} -\bar{x} } \leq \gamma^{2} \norm{x_{k} -\bar{x}} \leq \gamma^{2k+1}\norm{x_{0} -\bar{x}},
	\end{align*}
	where $\gamma:=\norm{\Pro_{ (D-\bar{x})\cap E^{\epsilon} }  \Pro_{ (C-\bar{x})\cap E^{\epsilon} }  }$ and
	$ E^{\epsilon} := \overline{(C-\bar{x})^{\ominus,\epsilon} +(D-\bar{x})^{\ominus,\epsilon} }$.
	
	Moreover, if $\hat{C}:=C-\bar{x}$ and $\hat{D}:=D-\bar{x}$ are convex cones, then $\gamma = c (\hat{C}, \hat{D})=c_{0} (\hat{C} \cap (\hat{C} \cap \hat{D}  )^{\ominus}, \hat{D} \cap (\hat{C} \cap \hat{D}  )^{\ominus}) =\norm{ \Pro_{\hat{D} \cap (\hat{C} \cap \hat{D}  )^{\ominus}}\Pro_{\hat{C} \cap (\hat{C} \cap \hat{D}  )^{\ominus}}}$.
\end{fact}

In order to study the cyclic projections algorithm for nonconvex sets, in \cite[Corollary~5.18]{LewisLukeMalick2009}, Lewis,  Luke and  Malick also show that under some conditions, if the initial point is sufficiently close to the solution, then the sequence of the cyclic projections algorithm converges with a   rate bounded below by a cosine of minimal angle between two cones, and bounded above by $1$. Again, this shows the importance of the study of the minimal angle between cones.

Notice that there are some papers (see, e.g., \cite{Obert1991}, \cite{SS2016I}, and \cite{Tenenhaus1988}) using the definition of \enquote{minimal angle} of two closed   convex cones with replacing the unit ball  by the unit sphere in   \cref{defn:Angles}. Although this \enquote{minimal angle} is consistent with our minimal angle defined in  \cref{defn:Angles} when it is in  $\left[0,\frac{\pi}{2}\right]$ (especially when the underlying two cones are linear subspaces),  its cosine is actually  in $\left[-1,1\right]$. But in this work we care mainly on using the cosine of minimal angle to describe the convergence rate of algorithms as applications presented above, and the convergence rate is always nonnegative. This is why we work only on the minimal angle defined in \cref{defn:Angles}.

In addition, although the definition of   \enquote{minimal angle} with cosine in $\left[-1,1\right]$ was used in \cite{Obert1991} and  \cite{Tenenhaus1988}, according to \cite[page~65]{Obert1991}, the minimal angle  in their problem in the real case is only located in $\left[0,\frac{\pi}{2}\right]$; and in view of \cite[Proposition~8]{Tenenhaus1988}, the optimal solution to the problem of searching for the minimal angle in  their canonical analysis associated with two convex polyhedral cones  can be obtained by a canonical analysis associated with two related linear subspaces. Note that applications of minimal angles  to ordinary differential equations and  to optimal multiple regression  were  presented in \cite{Obert1991} and \cite{Tenenhaus1988}, respectively.

We end   this work by summarizing our main results.  Recall that a pair of principal vectors is an optimal solution of a related minimal angle problem. We showed the existence of principal vectors of two nonempty convex sets (see \cref{prop:optimalpair:exist})  and provided  necessary conditions for a pair of vectors in $\mathcal{H} \times \mathcal{H}$ to be a pair of principal vectors of two nonempty convex sets  (see \cref{prop:optimalpair:neces}). In terms of the evaluation of minimal angle, we  presented equivalent expressions of the cosine of the minimal angle (see \cref{theorem:c0}) and characterized   the cosine of the minimal angle between two closed convex cones being strictly less than $1$ (see \cref{theorem:c0K1K2EQ}). We also specified sufficient conditions for the closedness of the sum of two nonempty closed convex cones in Hilbert spaces (see \cref{cor:K1K2closed}).
Moreover, we  proved that for two nonempty closed convex cones $K_{1}$ and $K_{2}$, under some conditions (e.g., the space is finite-dimensional), if  one of the two cones is not a linear subspace, then $c_{0}(K_{1},K_{2})<1$ and $c_{0}(K_{1}^{\ominus}, K_{2}^{\oplus}) <1$ cannot occur simultaneously (see \cref{coro:KominusMperpNeq0}).
Last but not least, some  counterexamples were constructed to confirm the tightness of our assumptions in related results (see \Cref{example:KMR3,exam:K1K2,exam:11NEQ} ).  
\section*{Acknowledgements}
The authors would like to thank the anonymous referees and the editor for their helpful comments.
HHB and XW were partially supported by NSERC Discovery Grants.

\addcontentsline{toc}{section}{References}

\bibliographystyle{abbrv}

\end{document}